\documentclass[12pt]{article}
\usepackage{amsmath,amssymb}
\newcommand{\eqdef}{\stackrel{\text{def}}{=}}
\newcommand{\n}{\nonumber}

\newcommand{\ignore}[1]{}
\numberwithin{equation}{section}
\newcommand{\Romannumeral}[1]{\uppercase\expandafter{\romannumeral#1}}
\newcommand{\I}{\text{\Romannumeral{1}}}
\newcommand{\II}{\text{\Romannumeral{2}}}
\newcommand{\ii}{(\text{\romannumeral2})}
\newcommand{\III}{\text{\Romannumeral{3}}}
\newcommand{\iii}{(\text{\romannumeral3})}
\newcommand{\IV}{\text{\Romannumeral{4}}}
\newcommand{\ai}{\text{I}}
\newcommand{\ait}{\text{II}}
\newcommand{\aitr}{\text{III}}

\newtheorem{prop}{\bf Proposition}[section]

\newcommand{\ns}{\hspace*{-7pt}}

\allowdisplaybreaks[4]

\begin{document}

\baselineskip=20pt

\newfont{\elevenmib}{cmmib10 scaled\magstep1}
\newcommand{\preprint}{
 \vspace*{-20mm}
    \begin{flushleft}
     \elevenmib Yukawa\, Institute\, Kyoto\\
   \end{flushleft}\vspace{-1.3cm}
   \begin{flushright}\normalsize \sf
     YITP-12-76\\
   \end{flushright}}
\newcommand{\Title}[1]{{\baselineskip=26pt
   \begin{center} \Large \bf #1 \\ \ \\ \end{center}}}
\newcommand{\Author}{\begin{center}
   \large \bf C.-L. Ho${}^{a,b}$, R. Sasaki${}^b$  and  K. Takemura${}^c$\end{center}}
\newcommand{\Address}{\begin{center}
  ${}^a$  Department of Physics, Tamkang University, Tamsui 25137, Taiwan\footnote{Permanent address}\\ 
${}^b$ Yukawa Institute for Theoretical Physics,
     Kyoto University,\\
      Kyoto 606-8502, Japan\\
  $^c$ Department of Mathematics, Faculty of Science and Technology,\\
      Chuo University, 
      1-13-27 Kasuga, Bunkyo-ku Tokyo 112-8551, Japan\\
      \end{center}}
\newcommand{\Accepted}[1]{\begin{center}
   {\large \sf #1}\\ \vspace{1mm}{\small \sf Accepted for Publication}
   \end{center}}

\preprint
\thispagestyle{empty}

\Title{Confluence of apparent singularities in multi-indexed orthogonal polynomials: the Jacobi case}

\Author

\Address
\vspace{1cm}

\begin{abstract}
The multi-indexed Jacobi polynomials  are the main part of the eigenfunctions of
exactly solvable quantum mechanical systems obtained by certain deformations of the P\"oschl-Teller potential
(Odake-Sasaki). By fine-tuning the parameter(s) of
the P\"oschl-Teller potential, we obtain several families of explicit and global solutions of
certain second order Fuchsian differential equations with an apparent singularity of characteristic exponent $-2$ and $-1$. 
They form orthogonal polynomials over $x\in(-1,1)$
with weight functions of the form $(1-x)^\alpha(1+x)^\beta/\{(ax+b)^4q(x)^2\}$, in which $q(x)$ is a polynomial in $x$.
\end{abstract}

\section{Introduction}
\label{sec:intro}

The theory of Fuchsian differential equations is essentially a local theory \cite{cod-lev,oshima}. Exactly solvable quantum mechanics, on the other hand, has recently provided infinitely many examples 
of global solutions, called  
the multi-indexed  Jacobi polynomials \cite{os25}, of second
order Fuchsian equations with $3+\ell$ ($\ell\ge1$) regular singularities. The first three regular singularities correspond to those for the Jacobi polynomials and the extra $\ell$ singularities are all apparent and their characteristic exponents are generically all $-1$. At this point a naive and deep question is `{\em if it is possible to construct global, {\em i.e.\/} rational solutions having higher ($\rho=-2,-3,\ldots$) characteristic exponents at apparent singularities?\/}' In a sense, this could be regarded as a problem of the confluence of apparent singularities.
In a previous paper \cite{st}, a related question was asked in connection with
the multi-indexed Laguerre polynomials \cite{os25,gomez3} and a 
partial positive answer was obtained. Namely, several examples
of second order differential equations with one irregular singularity at $\infty$ and two or more regular singularities at 
finite locations were constructed together with solutions having
a characteristic exponent $-2$ at  one of the regular singularities,
through fine-tuning of the parameter of the multi-indexed Laguerre polynomials.

In the present paper, we provide a partial positive answer to the
above question within the framework of the multi-indexed Jacobi
polynomials. As recapitulated in section two, the multi-indexed Jacobi polynomials are obtained as the main part of the eigenfunctions of a quantum system with the P\"oschl-Teller potential \cite{infhul,susyqm}, deformed by  Darboux-Crum
 transformations \cite{darb,crum,adler} in terms of the polynomial solutions called the {\em virtual state solutions\/}.
The ``confluence of apparent singularities" is achieved through fine-tuning of the two parameters $(g,h)$ of the Jacobi polynomial
to a curve $g=g(h)$ or $h=h(g)$ in such a way that the seed solutions have a cubic zero. In other words, these curves are the roots of the discriminants of the polynomial part of the candidate seed solutions, see \eqref{D-def}.
Here we show several examples of such Fuchsian equations (with a
$\rho=-2$) together with their global solutions.
We also restrict the parameter ranges such that the extra singularities do not appear in the physical region. 
Therefore the fine-tuned multi-indexed Jacobi polynomials form orthogonal polynomials over $(-1,1)$ with the weight function
of the form  $(1-x)^\alpha(1+x)^\beta/\{(ax+b)^4q(x)^2\}$, in which $q(x)$ is a polynomial in $x$. The denominator does not vanish
in $(-1,1)$.
Most of the explicit examples in section four have four regular singularities at $0$, $1$, $t$ and $\infty$, for which the Heun's equation is the canonical form:
\begin{equation}
\frac{d^2y}{dz^2} + \left( \frac{\gamma}{z}+\frac{\delta }{z-1}+\frac{\epsilon}{z-t}\right) \frac{dy}{dz} +\frac{\alpha \beta z -q}{z(z-1)(z-t)} y=0,
\label{Heun}
\end{equation}
with the condition  $\gamma +\delta +\epsilon =\alpha +\beta +1$
(Fuchs' relation). 
As is well known $(0,1-\gamma)$, $(0,1-\delta)$ and $(0,1-\epsilon)$ are the exponents
at $0$, $1$ and $t$, respectively and $(\alpha,\beta)$ are those
at $\infty$. The parameter $q$ is called the {\em accessory parameter\/}.
Heun's differential equation appears in several systems in physics \cite{Ron,SL}, 
in particular, its special cases appear 
as the differential equations for the exceptional Jacobi polynomials of type $X_1$, 
and some of their properties were studied in \cite{takemuraJPA}. 
By the change of the independent variable $\eta=1-2z$,  
the equations satisfied by the polynomials in all the examples of section four, \eqref{EoP12}, 
\eqref{EoP34}, \eqref{eopB1}, \eqref{eopB2}, \eqref{eopB3}, \eqref{eopB4} can be cast into the above canonical form \eqref{Heun}.

This paper is organised as follows. In section two, the essence of the multi-indexed Jacobi polynomials is recapitulated.
Starting with the general Schr\"odinger equation with  Darboux-Crum transformations in section \ref{sec:darb}, the details of the P\"oschl-Teller potential and the eigenpolynomials as well as the virtual state solutions are introduced in section
\ref{sec:dpt} and the explicit forms of the multi-indexed Jacobi polynomials are explained in section \ref{explicit}.
A r\'esum\'e of the method for the confluence of apparent singularities is given here.
The method to realise the higher degrees of apparent singularities is explained in section \ref{highsing}.
Section \ref{sec:fourgroups} is the main body of the paper. 
In section \ref{sec:four} eight distinct seed functions are chosen among the Wronskians of two and three virtual state solutions of lower degrees.
Their parameters are fine-tuned so that a cubic zero is realised.
They are divided into four groups, I, II, III and IV.  
The Groups I and II (III and IV) have the same form of potentials
\eqref{potI}, (respectively \eqref{potIII}) and thus the same form of differential 
equations \eqref{EoY12} and  \eqref{EoP12} 
(respectively \eqref{EoY34} and \eqref{EoP34}) when expressed in the shifted parameters $\bar{g}$ and $\bar{h}$ \eqref{ghbardef}. 
Within the same group, the potential and the differential equations of one member's can be obtained from the other's by a simple shift of the parameter, {\em e.g.}
$h\to h-2$, etc. In section \ref{sec:heunsys} the aspects of Heun's equation are discussed in some detail.
The limit to the multi-indexed Laguerre systems is discussed in section \ref{laglim}. Several other interesting cases are commented
on  in section \ref{others}.
The final section is for a summary and comments.
A summary table of the contents of section 4 and Appendix is provided here for 
reference purposes.
 Appendix \ref{datafour}
 provides the explicit forms of the fine-tuned multi-indexed polynomials, the differential equations for the rational as well as the polynomial solutions and the orthonormality relations for the Groups I, II, III and IV, together with the Laguerre limits of IV.

\section{Multi-indexed Jacobi polynomials}
\label{sec:miljac}
\setcounter{equation}{0}

Here we recapitulate the essence of the multi-indexed Jacobi polynomials as derived in \cite{os25}.
The main ingredients are (i) the Jacobi polynomials
as the main part of the eigenfunctions of the exactly solvable quantum system, the P\"oschl-Teller potential,
\ii\ the virtual state solutions,
\iii\ the  Darboux-Crum transformations.

\subsection{Darboux-Crum transformations}
\label{sec:darb}

We start with the Darboux-Crum transformation, which applies to any one-dimensional Schr\"odinger equation.
Let 
\begin{equation}
\mathcal{H}=-\frac{d^2}{dx^2}+U(x),
\label{schr}
\end{equation}
be a Hamiltonian and $U(x)$ be a potential defined in a complex domain.
To define the Darboux-Crum transformation, we pick up distinct seed solutions $\{\varphi_j(x), \tilde{\mathcal{E}}_j\}$ $(j=1,\ldots,M)$ of the original Schr\"odinger equation.
Namely,
\begin{equation}
 \mathcal{H}\varphi_j(x)
= \tilde{\mathcal{E}}_j\varphi_j(x),
\qquad  \tilde{\mathcal{E}}_j\in\mathbb{C},
\quad j=1,\ldots,M,
\label{scheq2}
\end{equation}
and the functions  $ \varphi_j(x)$  need not be square-integrable.
Let 
\begin{align}
&{\mathcal H}^{(M)}{\psi} (x)=\mathcal{E}{\psi} (x),
\label{Mschr} \\
&{\mathcal H}^{(M)}=-\frac{d^2}{dx^2}+{U}^{(M)}(x),\quad
{U}^{(M)}(x)\eqdef U(x)-2\frac{d^2\log\text{W}[\varphi_1,\ldots,\varphi_M](x)}{dx^2}
\label{Mwron}
\end{align}
be the $M$-th deformed Schr\"odinger equation, where $\text{W}[f_1,\ldots,f_N](x)$ is a Wronskian
\begin{equation}
\text{W}[f_1,\ldots,f_N](x)\eqdef 
\text{Det}\left(\frac{d^{k-1}f_j(x)}{dx^{k-1}}\right)_{1\le j,k\le N} .
\end{equation}
Then we have
\begin{prop} \label{prop:crumadler} {\rm \cite{crum,adler}}
Let $\psi(x) $ be a solution of the original Schr\"odinger equation 
\begin{equation}
\mathcal{H}\psi(x)=\mathcal{E}\psi(x) .
\label{scheq3}
\end{equation}
Then the function 
\begin{equation}
\psi^{(M)}(x)\eqdef\frac{\text{W}[\varphi_1,\ldots,\varphi_M,\psi](x)}
{\text{W}[\varphi_1,\ldots,\varphi_M](x)},
\label{psiM}
\end{equation}
satisfies the $M$-th deformed Schr\"odinger equation \eqref{Mschr}
with the same energy:
\begin{align}
&{\mathcal H}^{(M)}{\psi}^{(M)}(x)=\mathcal{E}{\psi}^{(M)}(x).
\end{align}
\end{prop}

In the $M$-th deformed potential \eqref{Mwron}, 
the zeros of the Wronskian ${\text{W}[\varphi_1,\ldots,\varphi_M](x)}$ are the extra singularities.

\subsection{P\"oschl-Teller potential}
\label{sec:dpt}

The Jacobi polynomials constitute the main part of the eigenfunctions of a quantum system
with the P\"oschl-Teller (PT) potential:
\begin{align}
\mathcal{H}=-\frac{d^2}{dx^2}+U(x;g,h), \quad 
U(x;g,h)=  {\displaystyle \frac{g(g-1)}{\sin^2x}+\frac{h(h-1)}{\cos^2x}-(g+h)^2},
\label{pot}\\
g>0,\quad h>0,\quad x\in \left(0, \frac{\pi}{2}\right).
\end{align}
It has regular singularities at $x=0 $ and  $x=\pi /2 $ modulo $\pi$,
and the characteristic exponents at $x=0 $ (resp. $x=\pi /2 $)  are
$\rho=g,1-g$ (resp. $\rho=h,1-h$).
The Hamiltonian is self-adjoint and the eigenfunctions and the eigenenergies are given by
\begin{align}
\phi_n(x;g, h)&=(\sin x)^g (\cos x)^h P^{(g - 1/2, h - 1/2)}_n(\eta(x)),\quad 
\eta(x)\eqdef\cos(2x),\\
 \mathcal{E}_n(g,h)&= 4 n(n+g+h),\qquad \qquad
  n=0,1,2, \dots,
\end{align}
where   $P_n^{(\alpha,\beta )}(\eta )$ is the Jacobi polynomial 
in the variable $\eta$ defined by
\begin{equation}
P _n^{(\alpha , \beta )}(\eta )= \frac{(\alpha +1)_n}{n!} \sum _{k=0}^n \frac{(-n)_k(n+\alpha +\beta +1)_k}{k! (\alpha +1) _k} \left( \frac{1-\eta }{2} \right)^k .
\label{jacobi}
\end{equation}
The eigenpolynomial solutions form a complete set of orthogonal polynomials over $(0,\pi /2)$ under a suitable restriction of the parameters $g,h$:
\begin{align} 
   (\phi_n,\phi_m)&\eqdef \int_0^{\pi /2} \! \phi_n(x;g,h)\phi_m(x;g,h)dx =h_n(g,h)\delta_{n\,m},
   \label{intformori}\\
   &\qquad \quad 
    h_n(g,h)\eqdef \frac{\Gamma(n+g+\frac12)\Gamma(n+h+\frac12)}
  {2n!(2n+g+h)\Gamma(n+g+h)}.
    \label{hform}
\end{align}
Another ingredient is the {\em virtual state solutions\/}, which are three types of seed polynomial solutions of the Schr\"odinger equation \eqref{pot} indexed by $\text{v}\in\mathbb{Z}_{>0}$, which have all real energies:
\begin{align}
    \tilde{\phi}_\text{v}^\ai(x;g,h)&\eqdef 
(\sin x)^g (\cos x)^{1 - h} P^{(g - 1/2, 1/2-h)}_\text{v}(\eta(x)),~~~\n\\
&~~~~~~\tilde{\mathcal{E}}_\text{v}^\ai(g,h)=-4 (g + \text{v} + 1/2) (h   -\text{v}-1/2);
 \label{seed1}\\
  \tilde{\phi}_\text{v}^\ait(x;g,h)&\eqdef 
(\sin x)^{1-g} (\cos x)^h P^{(1/2-g, h-1/2)}_\text{v}(\eta(x)),~~~\n\\
&~~~~~~\tilde{\mathcal{E}}_\text{v}^\ait(g,h)= -4 (g -\text{v}- 1/2) (h +  \text{v}+ 1/2);
    \label{seed2}\\
 \tilde{\phi}_\text{v}^\aitr (x;g,h)&\eqdef  
 (\sin x)^{1-g} (\cos x)^{1 - h}  P^{(1/2-g, 1/2-h)}_\text{v}(\eta(x)),~~~\n\\
 &~~~~~\tilde{\mathcal{E}}_\text{v}^\aitr(g,h)=  -4 (\text{v}+ 1) (g + h - 1 -\text{v} )=\mathcal{E}_{-(\text{v}+1)}(g,h).
 \label{seed3}
\end{align}
Note that the seed solutions of type $\ai$ and $\ait$ are used to construct generic multi-indexed Jacobi polynomials in \cite{os25}. 
The third type $\aitr$ is now employed to produce a wider class of polynomials, to be called tentatively {\em non-generic multi-indexed Jacobi polynomials\/} (for the non-generic
multi-indexed Laguerre polynomials, see \cite{st}). Because of the type III virtual solutions, the absence of singularities in $\eta\in(-1,1)$ is no longer guaranteed.
One salient feature of the non-generic multi-indexed polynomials is the existence of the additional modes at $n<0$.
When $\tilde{\phi}^\aitr_\text{v}$ is used as a seed solution, an additional mode at $n=-(\text{v}+1)$ is created.
The degree of the polynomial $\tilde{\phi } ^{\ai}_{\rm v} (x) $ (resp. $\tilde{\phi } ^{\ait}_{\rm v} (x) $, $\tilde{\phi } ^{\rm III}_{\rm v} (x) $) is equal to $\rm v$ under the condition 
\begin{equation}
g-h+k \neq 0 \ (\text{resp.}\ g-h-k \neq 0, g+h+1-k \neq 0)\ \text{for}\ k={\rm v}  +1, {\rm v} +2, \dots ,2{\rm v}.
\end{equation}
In the rest of this section, we assume that the parameters are {\em general} and satisfy these conditions, {\em unless otherwise stated\/}.
In section \ref{highsing} and after, fine tuning of the parameters will be used extensively and the above conditions will not be assumed any longer.

%
%
\subsection{Explicit forms of the multi-indexed Jacobi polynomials}
\label{explicit}

We define a wider class of  multi-indexed Jacobi polynomials by applying Darboux-Crum transformation with  {\em three\/} types of the seed solutions to the quantum system with the PT potential in \eqref{pot}. Let $\tilde{\phi } ^{\text{t}_1}_{\text{v}_1},\ldots,\tilde{\phi } ^{\text{t}_\mathcal{N}}_{\text{v}_\mathcal{N}} $ $(\text{t}_1, \dots \text{t}_\mathcal{N} \in \{ \ai , \ait, \aitr \} , \;  \text{v}_1, \dots \text{v}_\mathcal{N} \in {\mathbb Z }_{>0})$ be distinct seed solutions. 
Let us denote the number of type I seed solutions in the set
$\{\tilde{\phi } ^{\text{t}_1}_{\text{v}_1},\ldots,\tilde{\phi } ^{\text{t}_\mathcal{N}}_{\text{v}_\mathcal{N}}\}$  by $M\ge0$,  that of type II  by $N\ge0$ and that of type III  by $L\ge0$.
Then $\mathcal{N}\eqdef M+N+L$ is the total number of the used seed solutions. 
It follows from Proposition \ref{prop:crumadler} that the function 
\begin{equation}
\phi^{(\mathcal{N})}_n(x)\eqdef
\frac{\text{W}[\tilde{\phi } ^{\text{t}_1}_{\text{v}_1},\ldots,\tilde{\phi } ^{\text{t}_\mathcal{N}}_{\text{v}_\mathcal{N}},\phi_n](x)}{\text{W}[\tilde{\phi } ^{\text{t}_1}_{\text{v}_1},\ldots,\tilde{\phi } ^{\text{t}_\mathcal{N}}_{\text{v}_\mathcal{N}}](x)}
\label{eq:phiMn}
\end{equation}
 is an eigenfunction of the deformed Hamiltonian 
\begin{align}
&{\mathcal H}^{(\mathcal{N})}=-\frac{d^2}{dx^2}+ U(x;g,h)-2\frac{d^2\log\text{W}[\tilde{\phi } ^{\text{t}_1}_{\text{v}_1},\ldots,\tilde{\phi } ^{\text{t}_\mathcal{N}}_{\text{v}_\mathcal{N}}](x)}{dx^2}
\end{align}
with the same eigenvalue $\mathcal{E}_n=4n(n+g+h)$, provided that
the deformed potential is non-singular. 
Let us express the numerator and the denominator of $\phi^{(\mathcal{N})}_n(x)$ as
\begin{align}
\text{W}[\tilde{\phi } ^{\text{t}_1}_{\text{v}_1},\ldots,\tilde{\phi } ^{\text{t}_\mathcal{N}}_{\text{v}_\mathcal{N}},\phi_n](x)
&\eqdef (1-\eta)^\Lambda(1+\eta)^\Sigma \mathcal{Q}_n(\eta),\quad
\mathcal{Q}_n(\pm1)\neq0,
\label{Q-def}\\
\text{W}[\tilde{\phi } ^{\text{t}_1}_{\text{v}_1},\ldots,\tilde{\phi } ^{\text{t}_\mathcal{N}}_{\text{v}_\mathcal{N}}](x)
&\eqdef (1-\eta)^{\Lambda'}(1+\eta)^{\Sigma'} \mathcal{D}(\eta),\quad
\mathcal{D}(\pm1)\neq0,
\label{D-def}
\end{align}
in which $\mathcal{Q}_n(\eta)$ and $\mathcal{D}(\eta)$ are polynomials in $\eta$.
\begin{prop} \label{prop:multip}
The exponents are given by
\begin{align}
\Lambda&\eqdef \frac14(J+1)(2g+J),\quad 
\Sigma\eqdef\frac14(K+1)(2h+K),\\
\Lambda'&\eqdef \frac14 J(2g+J-1),\quad \ \  
\Sigma'\eqdef\frac14 K(2h+K-1),\\
\text{with}\quad J&\eqdef M-N-L,\qquad\quad  K\eqdef -M+N-L.
\end{align}
The deformed potential is given by
\begin{align}
U^{(\mathcal{N})}(x)&\eqdef \frac{2\bar{g}(\bar{g}-1)}{1-\eta}+\frac{2\bar{h}(\bar{h}-1)}{1+\eta}-(\bar{g}+\bar{h}+2L)^2\n\\
&\qquad \qquad \qquad\qquad 
+8\left[ (1-\eta ^2) \frac{\mathcal{D}'(\eta ) ^2-\mathcal{D}(\eta) \mathcal{D}''(\eta) }{\mathcal{D}(\eta)^2} +\eta \frac{\mathcal{D}'(\eta)}{\mathcal{D}(\eta )} \right],
\label{defDPT}\\
&\bar{g}\eqdef g+M+N-L,\quad
\bar{h}\eqdef h-M+N-L,\quad \bar{g}+\bar{h}+2L=g+h,
\label{ghbardef}
\end{align}
and the eigenfunctions
have a simple form
\begin{align}
\phi^{(\mathcal{N})}_n(x)&=(1-\eta)^{\bar{g}/2}(1+\eta)^{\bar{h}/2}
\frac{\mathcal{Q}_n(\eta)}{\mathcal{D}(\eta)},\quad \text{\rm deg}(\mathcal{Q}_n)=\ell+n+L,\quad 
\text{\rm deg}(\mathcal{D})=\ell,
\label{poldef}\\
\ell&\eqdef \sum_{j=1}^{M+N+L}\text{\rm v}_j-\frac{M(M-1)}2-\frac{N(N-1)}2-\frac{L(L-1)}2+MN\ge0.
\label{deg}
\end{align}
For general values of $g$ and $h$, the two polynomials $\mathcal{Q}_n(\eta)$ and $\mathcal{D}(\eta)$ are relatively prime and the polynomials $\{\mathcal{Q}_n(\eta)\equiv \mathcal{Q}_n(\eta;g,h)\}$ are called the 
multi-indexed Jacobi polynomials.
\end{prop}
With a slight change of notation, these formulas are generalisation of those in \cite{os25} with the inclusion of the type III seed solutions.
For example, the formulas \eqref{ghbardef}, \eqref{poldef} and \eqref{deg} correspond to (26), (25) and (30) of \cite{os25}.

The regularity of the solutions and/or the positivity of the resulting weight 
functions etc. must be verified in each specific case for the chosen parameter values. This is in good contrast with the generic multi-indexed Laguerre and Jacobi polynomials in \cite{os25}, in which the conditions for the parameters are given (eqs.(23) and (24) of \cite{os25}).

When there is no singularity in $\eta\in(-1,1)$,
the norm of the $\mathcal{N}$-th transformed solution $\phi_n^{(\mathcal{N})}$  of the $n$-th the eigenfunction $\phi_n(x)$ 
is related to that of the original $\phi_n(x)$:
\begin{align}
&(\phi^{(\mathcal{N})}_n,\phi^{(\mathcal{N})}_n)=
\prod_{j=1}^M(\mathcal{E}_n- \tilde{\mathcal{E}}^{\I}_{\text{v}_j})\cdot
\prod_{j=1}^N(\mathcal{E}_n- \tilde{\mathcal{E}}^{\II}_{\text{v}_j})
\cdot\prod_{j=1}^L(\mathcal{E}_n- \tilde{\mathcal{E}}^{\III}_{\text{v}_j})\cdot(\phi_n,\phi_n),
\label{intformM}
\end{align}
where  $\tilde{\mathcal{E}}_j$ is the energy of the $j$-th seed solution.

\bigskip
Here is a r\'esum\'e of the method for the confluence of apparent singularities.
Based on the general Darboux transformation formulas \eqref{Mschr}--\eqref{psiM},
the multi-indexed Jacobi polynomials $\mathcal{Q}_n(\eta)$  \eqref{poldef} are derived
from the eigenfunctions of the P\"oschl-Teller potential \eqref{pot} by using the three
types of virtual state wavefunctions \eqref{seed1}--\eqref{seed3}.
A high degree of apparent singularity is realised when the denominator polynomial 
$\mathcal{D}(\eta)$ \eqref{poldef}  determined by \eqref{D-def} has a multiple zero \eqref{mzero}.
As explained in the next section, a multiple zero is achieved through fine-tuning  of
the parameters of the original Jacobi polynomials \eqref{jacobi} at the zeros of the discriminant of the denominator polynomial 
$\mathcal{D}(\eta)$ \eqref{poldef}.

%
%
\section{Higher degree of apparent singularities}
\label{highsing}
\setcounter{equation}{0}
We are going to investigate the relationship between the zero of the Wronskian and the singularity of the deformed Schr\"odinger equation.

Let us start with a general Schr\"odinger equation \eqref{schr} with
an unspecified potential $U(x)$.
Suppose  $x_0$, a regular point of the original potential $U(x)$, is a zero of the Wronskian  
$\text{W}[\varphi_1,\ldots,\varphi_\mathcal{N}](x)$ \eqref{Mwron}  with 
multiplicity $m$, {\em i.e.} 
\begin{equation}
\text{W}[\varphi_1,\ldots,\varphi_\mathcal{N}](x) =c_0(x-x_0)^m+O(x-x_0)^{m+1},\quad m\in\mathbb{Z}_{>0}.
\label{mzero}
\end{equation}
Then it is a singularity of the potential of the deformed Hamiltonian ${\mathcal H}^{(\mathcal{N})}$ such that
\begin{equation}
x\approx x_0,\qquad {U}^{(\mathcal{N})}(x)=U(x_0)
+\frac{2m}{(x-x_0)^2}+\text{regular terms}.
\end{equation}
The point $x_0$ is an apparent singularity of the deformed Schr\"odinger equation, 
because the monodromy is preserved by  Darboux transformations.
The characteristic exponents at the apparent singularity $x_0$ are determined by
\begin{equation}
\rho(\rho-1)-2m=0,
\end{equation}
and the difference of the exponents is an integer because of the apparency of the singularity.
Combining it with that $m \in {\mathbb Z} _{>0}$, the admitted values are  
\begin{align} 
   m&=l (l +1)/2 ,\quad \rho=-l,\quad l +1,
   \quad (l \in {\mathbb Z}_{>0}).
\end{align}
Hence the degrees of zeros of the Wronskian are restricted to 
$m=1, 3, 6, 10, \dots $.
In the case $m=1$,  the exponents are $-1$ and $2$, and the higher 
exponents of the singularity are realised only in the cases $m=3,6,10, \dots $.

Next we investigate the condition that the Schr\"odinger equation of 
the deformed P\"oschl-Teller potential \eqref{defDPT} has 
an extra singularity in $\eta $ with higher exponents 
on a curve $g=g(h)$ (or equivalently, $h=h(g)$) in the two dimensional
parameter plane $(g,h)$.
This means that the polynomial $\mathcal{D}(\eta; g(h), h)$ has a multiple zero except for $\eta =\pm 1$.
Hence a necessary condition for the curve $g=g(h)$ is that the discriminant of the polynomial $\mathcal{D}(\eta;g(h),h)$ in the variable $\eta $ vanishes.
It is not sufficient because the multiple zeros of $\mathcal{D}(\eta ; g(h), h)  $ may be realised at $\eta =\pm 1$ or the  coefficient of the leading term of  $\mathcal{D}(\eta ;g(h), h)$ may be zero.
If the polynomial $\mathcal{D}(\eta ;g(h) , h)$ has a multiple zero at $\eta = \eta _0(h)$, then the degree of the zeros are restricted to $m=3, 6, 10, \dots $.
By fine-tuning the parameters on the curve $g=g(h)$,  the confluence of some  apparent singularities is achieved. 
The fine-tuning has another effect. At the multiple zeros,  the 
numerator polynomial $\mathcal{Q}(\eta)$ and the denominator polynomial
$\mathcal{D}(\eta)$ are {\em no longer relatively prime.\/}
While the denominator $\mathcal{D}(\eta)$ has an $m=l(l+1)/2$ zero, the numerator polynomial $\mathcal{Q}(\eta)$
has an $l(l+1)/2-l=l(l-1)/2$ zero:
\begin{align} 
  \mathcal{D}(\eta;g(h),h)&=(\eta-\eta_0(h))^{l(l-1)/2}w(\eta;g(h),h),
  \label{wdef} \\
 w(\eta;g(h),h)&=(\eta-\eta_0(h))^l\times q(\eta),\quad \text{deg}(q)= \ell-l,\\
  \mathcal{Q}_n(\eta;g(h),h)&=(\eta-\eta_0(h))^{l(l-1)/2}\mathcal{P}_n(\eta;g(h),h),
   \quad \mathcal{P}_n(\eta_0(h);g(h),h)\neq0.
\end{align}
The function $\phi^{(\mathcal{N})}_n(x) = (1-\eta )^{\bar{g}(h)/2 }(1+\eta )^{\bar{h}(h)/2} \mathcal{P}_n(\eta;g(h),h)/w(\eta;g(h),h)$ defined by \eqref{poldef} satisfies the deformed  Schr\"odinger equation and the exponents are $-l $ and $l +1$.
If $l =2$,  the denominator Wronskian has a triple zero, and the degree one polynomial is cancelled  by that of the numerator.
It should be stressed that the normalisation of the polynomials $\{\mathcal{P}_n\}$  is uniquely determined by that of $ w(\eta;g(h),h)$.

Hereafter, for simplicity of presentation, we simply write 
$\mathcal{P}_n(\eta;g,h)$, etc. when the curve $g=g(h)$ (or $h=h(g)$) is pre-specified.

 %
%

\section{Four groups of new Jacobi  systems}
\label{sec:fourgroups}
\setcounter{equation}{0}

Now we present several concrete examples of interesting systems
exhibiting apparent singularities with the characteristic exponent $\rho=-2$. 
We shall be interested in new quantum systems generated by the  Darboux transformation using the following  
eight sets of  seed functions:
\begin{alignat*}{4} 
& [ \tilde{\phi}_2^\ai,\tilde{\phi}_1^\aitr] ;&\hspace{10mm} &[ \tilde{\phi}_1^\ai,\tilde{\phi}_2^\aitr];&\hspace{10mm}
&  [ \tilde{\phi}_1^\ai, \tilde{\phi}_2^\ai,\tilde{\phi}_1^\aitr] ; 
&\hspace{10mm} &[ \tilde{\phi}_1^\ai, \tilde{\phi}_1^\aitr,\tilde{\phi}_2^\aitr] ;\\
& [ \tilde{\phi}_2^\ait,\tilde{\phi}_1^\aitr] ; &\hspace{10mm} &[ \tilde{\phi}_1^\ait,\tilde{\phi}_2^\aitr] ; &\hspace{10mm}
&[ \tilde{\phi}_1^\ait, \tilde{\phi}_2^\ait,\tilde{\phi}_1^\aitr]; 
&\hspace{10mm}  &[ \tilde{\phi}_1^\ait, \tilde{\phi}_1^\aitr,\tilde{\phi}_2^\aitr].\n
\end{alignat*} 
The four sets in the second row are obtained by replacing in the first row the type I seed functions by the type II ones.
Systems generated by other choices of the set of seed functions can be studied accordingly by the method described in the main text.

\subsection{The four groups of Jacobi systems}
\label{sec:four}
We look for the situation where the polynomial part of the 
Wronskian of the set of seed functions, $\mathcal{D}(\eta)$ \eqref{D-def},  contains a cubic zero outside of the physical domain 
$x\in (0,\pi/2)$ (or $\eta\in(-1,1)$).
It turns out that such a situation exists for these eight cases above when $g=g(h)$ is a linear fractional function.  
Note that it is equally well to consider $h=h(g)$ as a function of $g$ in each case. But taking $g$ as a function of $h$
allows us to consider the $h\to \pm\infty$ limit, which will take the P\"oschl-Teller potential to the radial oscillator 
 considered in \cite{st}.

It is discovered out  that these eight cases fall into four groups of two members each,  such that  one member of a 
group can be obtained from the other by a shift in the value of $h$, which is determined by the universal formula for $\bar{g}$ and 
$\bar{h}$, \eqref{ghbardef}  derived in section \ref{explicit}.   
This is because, the main part of the denominator polynomial $\mathcal{D}(\eta;g(h),h)$, which has a cubic zero, is the same for the two members when expressed in terms of $\bar{h}$.
Moreover, the polynomial up to a multiplicative constant is identical
for the Groups I and II and it has the same form for the Groups III and IV, with only the value of a parameter $\gamma$ is different.
This means that the deformed potential $U^{(\mathcal{N})}(x)$ is identical for the Groups I and II, and it has the same form in the Groups III and IV, {\em when expressed in terms of $\bar{g}$ and $\bar{h}$}.

\begin{enumerate}
\item[$\bullet$] 
 Group I: (a)  $[\tilde{\phi}_2^\ai,\tilde{\phi}_1^\aitr]$,
 \quad \quad $g=\frac{3(h-3)}{4h-9},\quad \ \bar{g}=g,\qquad  \ \bar{h}=h-2$; \\
 Group I: (b) $[\tilde{\phi}_1^\ait, \tilde{\phi}_2^\ait,\tilde{\phi}_1^\aitr]$, $g=\frac{3(5h+3)}{4h+3}, \quad \bar{g}=g-3, \quad \bar{h}=h+1$.\\
 For both cases,  $\bar{g}=\frac{3(\bar{h}-1)}{4\bar{h}-1}$, $\bar{g}+\bar{h}+2=g+h$, ($L=1$).
\item[$\bullet$] 
 Group II: (a) $[\tilde{\phi}_1^\ai,\tilde{\phi}_2^\aitr]$,
 \quad  $g=\frac{h}{4h-9},\quad  \bar{g}=g, \quad  \bar{h}=h-2$;\\
 Group II:  (b) $[\tilde{\phi}_1^\ait,\tilde{\phi}_2^\aitr]$,  \quad $g=\frac{9h}{4h-1},\quad  \bar{g}=g-2, \quad  \bar{h}=h$. \\
  For both cases,  $\bar{g}=\frac{\bar{h}+2}{4\bar{h}-1}$, $\bar{g}+\bar{h}+2=g+h$, ($L=1$).
\begin{align} 
 \text{For Groups I \& II:}\qquad &\mathcal{D}_{\I,\II}(\eta;g,h)=(2\bar{h}\eta-2\eta-2\bar{h}-1)^3\times\text{constant},
 \label{D1}\\
 &U^{(\mathcal{N})}_{\I,\II}(x)=\frac{2\bar{g}(\bar{g}-1)}{1-\eta}+\frac{2\bar{h}(\bar{h}-1)}{1+\eta}-(\bar{g}+\bar{h}+2)^2\qquad \qquad
 \qquad \n\\[2pt]
& \qquad -8\left[\frac{3(2\bar{h}+1)}{2\bar{h}\eta-2\eta-2\bar{h}-1}
+\frac{9(4\bar{h}-1)}{(2\bar{h}\eta-2\eta-2\bar{h}-1)^2}\right],
\label{potI}
\end{align}
\item[$\bullet$]
 Group III:  (a) $[\tilde{\phi}_2^\ait,\tilde{\phi}_1^\aitr]$, 
 \qquad \ $g=\frac{9(h-1)}{4h-3}, \quad  \bar{g}=g-2,\quad  \bar{h}=h$; \\
 Group III:  (b) $[ \tilde{\phi}_1^\ai, \tilde{\phi}_2^\ai,\tilde{\phi}_1^\aitr]$,
 \quad $g=-\frac{3(h-3)}{4h-15},\quad  \bar{g}=g+1,\quad  \bar{h}=h-3$. \\
  For both cases,\quad  $\bar{g}=\frac{\bar{h}-3}{4\bar{h}-3}$,\quad 
  \ $\bar{g}+\bar{h}+2=g+h$, ($L=1$).
 \item[$\bullet$]
 Group IV:  (a)  $[ \tilde{\phi}_1^\ait, \tilde{\phi}_1^\aitr,\tilde{\phi}_2^\aitr]$,
 \quad $g=\frac{3(5h-8)}{4h-7},\quad \bar{g}=g-3,\quad \bar{h}=h-1$; \\
 Group IV: (b) $[ \tilde{\phi}_1^\ai, \tilde{\phi}_1^\aitr,\tilde{\phi}_2^\aitr]$,
 \quad $g=\frac{7h-24}{4h-15}, \quad \bar{g}=g-1,\quad  \bar{h}=h-3$. \\
 For both cases, \quad $\bar{g}=\frac{3\bar{h}}{4\bar{h}-3}$, \quad
 $\bar{g}+\bar{h}+4=g+h$, ($L=2$).
\begin{align} 
 \text{For Groups III \& IV:}\qquad &\mathcal{D}_{\III,\IV}(\eta;g,h)
 =(2\bar{h}\eta-2\bar{h}+3)^3\times\text{constant},\\
 &U^{(\mathcal{N})}_{\III,\IV}(x)=\frac{2\bar{g}(\bar{g}-1)}{1-\eta}+\frac{2\bar{h}(\bar{h}-1)}{1+\eta}-(\bar{g}+\bar{h}+\gamma)^2\qquad \qquad
 \quad \n\\[2pt]
& \qquad\qquad -8\left[\frac{3(2\bar{h}-3)}{2\bar{h}\eta-2\bar{h}+3}
-\frac{9(4\bar{h}-3)}{(2\bar{h}\eta-2\bar{h}+3)^2}\right],
\label{potIII}\\[2pt]
&~~~~~  \gamma=2~ (\text{for \   Group III}), ~4 ~(\text{for \   Group IV}).\label{U2}
\end{align}  
\end{enumerate}
The normalisation of the (fine-tuned) multi-indexed polynomials $\{\mathcal{P}_n(\eta)\}$ is fixed by 
\begin{align} 
\phi^{(\mathcal{N})}_n(x) &= (1-\eta )^{\bar{g}(h)/2 }(1+\eta )^{\bar{h}(h)/2} \mathcal{P}_n(\eta;g(h),h)/w(\eta;g(h),h),
\label{yndef}\\
w(\eta)&= w_1(\eta;\bar{g}, \bar{h})\eqdef (2\bar{h}\eta-2\eta-2\bar{h}-1 )^2,\quad \text{for Groups I \& II},\n\\
w(\eta)&= w_2(\eta;\bar{g}, \bar{h})\eqdef (2\bar{h}\eta-2\bar{h}+3)^2,\qquad \quad \text{for Groups III \& IV}.\n
\end{align} 
 For these to be non-singular in the physical domain $\eta\in (-1,1)$, one must restrict 
\begin{align} 
& \bar{h}>\frac14\quad {\rm and} \quad   \bar{h}\neq 1,  : \text{I,\  II},\qquad
 \bar{h}<\frac34\quad {\rm and} \quad   \bar{h}\neq 0 : \text{III,\  IV}.
\label{h-range}
\end{align} 
The conditions $ \bar{h}\neq 1, 0$ are also needed, 
otherwise  the
denominator in each case will be a constant, which is not the case we want to consider here.
The square integrability of the eigenfunctions requires
\begin{equation}
\bar{g}>-\frac12,\quad \bar{h}>-\frac12.
\label{squarecond}
\end{equation}
In Appendix  we list the possible range of $h$  for each specific case.

  We note here that  $\bar{g}_{\text{I}}=1-\bar{g}_{\text{II}}$ and $\bar{g}_{\text{III}}=1-\bar{g}_{\text{IV}}$.  This explains why the potentials in $\eta$ 
  are the same for Groups I and II, and for III and IV.
By equating $\bar h$ for the two cases in the same group, one obtains a transformation  between the two systems.   
Hence by setting $h\to h+3$ in the potential, the eigenvalues, the eigenfunctions etc in Case Ia (with $g$ expressed in 
terms of $h$), one obtains the corresponding formulas for Case Ib. Similarly, setting $h\to h+2$, $h\to  h-3$ and
 $h\to h-2$ map Cases IIa to IIb, IIIa to IIIb and IVa to IVb, respectively. These are the required shiftings mentioned above. 
 
 %
 %
 
Here are the Fuchsian differential equations for the rational functions
$y_n(\eta)\eqdef\mathcal{P}_n(\eta)/w(\eta)$, $n=0,1,\ldots,$ with a $\rho=-2$ apparent singularity.  They are obtained by rewriting the Schr\"odinger equation with the potential $U^{(\mathcal{N})}(x)$ \eqref{potI}, \eqref{potIII}:
\begin{align} 
&  (1-\eta^2) y^{\prime\prime}_n(\eta) 
+\left(\bar{h}-\bar{g} -(\bar{g} +\bar{h} +1)\eta\right) y^\prime_n(\eta)  +\left\{\phantom{\text{\Large I}}\!\!
(n+1)(n+\bar{g}+\bar{h}+1)\right. \n\\
&\left. \qquad \qquad \qquad \quad +\frac{6(2\bar{h}+1)}{2\bar{h}\eta-2\eta-2\bar{h}-1}+\frac{18(4\bar{h}-1)}{(2\bar{h}\eta-2\eta-2\bar{h}-1)^2}\right\}y_n(\eta)=0,
\label{EoY12}
\end{align} 
for Groups I and II,
and
\begin{align} 
&  (1-\eta^2) y^{\prime\prime}_n(\eta) 
+\left(\bar{h}-\bar{g} -(\bar{g} +\bar{h} +1)\eta\right) y^\prime_n(\eta)  +\left\{\phantom{\text{\Large I}}\!\!
(n+\frac{\gamma}2)(n+\bar{g}+\bar{h}+\frac{\gamma}2)\right. \n\\
&\left. \qquad \qquad \qquad \qquad \qquad \quad \ +\frac{6(2\bar{h}-3)}{2\bar{h}\eta-2\bar{h}+3}-\frac{18(4\bar{h}-3)}{(2\bar{h}\eta-2\bar{h}+3)^2}\right\}y_n(\eta)=0,
\label{EoY34}
\end{align}
for  Groups III and IV.   In these equations, $\bar{g}=\bar{g}(\bar{h})$ should be replaced by its respective value for each group as given before.
 These are second order differential equations  with regular singularities
 at $\eta=\pm 1$, $\infty$ and  
 $\eta=(2\bar{h}+1)/2(\bar{h}-1)$ 
(for Groups I and II), or $1-3/2\bar{h}$ (for Groups III, IV) with the 
 conditions (\ref{h-range}).

The polynomials $\mathcal{P}_n$ satisfy   
\begin{align} 
    & \left\{1+2n+n^2-2 \bar{h}^2 (3 + n) (\eta-1) - 6 \eta + 
     4n\eta + 2 n^2\eta \right.\n\\
&\quad + 
     \bar{h} \left[-1+n(5-2\eta)- 2n^2(\eta-1)+ 12 \eta\right] + 
     \bar{g} \left[5+n+6\eta+2n\eta \right. \n\\
&\left. \left. \quad  
   -2\bar{h}(1+ n(\eta-1)+3 \eta)\right]\right\}\mathcal{P}_n(\eta) + \left\{-8+2\bar{h}^2(\eta-1)^2-\eta+6\eta^2\right. \n\\
&\left. \qquad \quad + 
     \bar{g}(-1+ 2\bar{h}(\eta-1)-2\eta)(\eta+1) 
     -\bar{h}(\eta-1)(8\eta+9)\right\}\mathcal{P}_n^\prime(\eta)
     \n\\
&\qquad \qquad  \qquad \qquad  \qquad 
     + (-1+ 2\bar{h}(\eta-1) -2 \eta)(\eta^2-1) 
    \mathcal{P}^{\prime\prime}_n(\eta)=0,
    \label{EoP12}
\end{align}
 for Groups I and II, and
\begin{align} 
-\left\{3(\gamma+2n)^2 + 
      2\bar{h}\left[3(\gamma+2n)-16\eta+(\gamma+2n)^2  (\eta-1) \right]+ 
      4\bar{h}^2(4+\gamma+2 n)(\eta-1)\right.\n\\
 \left.  + 
      2 \bar{g} \left[3 (\gamma + 2 n)
 + 2 \bar{h} (4(\eta +1) +(\gamma+2 n)(\eta-1)
        )\right]\right\} \mathcal{P}_n(\eta) \n\\
            + 
  4 \left\{(2\bar{h}^2(\eta-1)^2 + 3\eta + 
        \bar{g}\left[3 + 2 \bar{h} (\eta-1)\right] ( \eta+1) - 
        \bar{h}(6\eta+5)(\eta-1)\right\} \mathcal{P}_n^\prime(\eta)\n\\
        + \left[3+2\bar{h} (\eta-1)\right] (\eta^2-1) 
\mathcal{P}_n^{\prime\prime}(\eta)=0,
        \label{EoP34}
\end{align}
 for Groups III and IV.
       
 Appendix  \ref{datafour}
 presents the essential data for these new polynomials for one case from each 
group.  The other member can be obtained by the shifting mentioned in the last paragraph.
Examples IIb and IIIa were discussed in \cite{gomez4} from a
different point of view.

\subsection{Heun's systems}
\label{sec:heunsys}
All the examples in the previous subsection have four regular singularities,
three at finite locations and one at $\infty$.
Thus they belong to the Heun's systems with one apparent singularity of characteristic exponent $\rho=-2$. Let us discuss the Heun's properties by taking
the example, Group I (a) with $g=3(h-3)/(4h-9)$. The $\mathcal{D}(\eta)$ 
polynomial \eqref{D1}  in terms of $h=\bar{h}+2$ reads
\begin{equation*}
\mathcal{D}(\eta;g,h)=( 2(h-3)\eta -2h+3 )^3\times\text{constant}.
\end{equation*}
By the change of the independent variable $\eta=1-2z$, the equation for the
two-indexed Jacobi polynomial $f_n=\mathcal{P}_n (\eta ;3(h-3)/(4h-9),h)$ is rewritten into the Heun's differential equation
\eqref{Heun}:
\begin{align}
 & \frac{d^2 f_n}{dz ^2}+\left( \frac{g+1/2}{z}+\frac{h-3/2}{z-1 }-\frac{4}{z+ \frac{3}{ 4(h-3)}} \right) \frac{df_n}{dz }  
 \n\\[4pt]
& + \frac{(-n-3) (n+g+h-3)z -\frac{3(4 h-9 )n^2+3 (4 h^2-6 h-9 )n+2 h (10 h-27)}{ 4(h-3) (4 h-9)}}{ z(z-1) (z+ \frac{3}{ 4(h-3)} )} f_n=0,
\quad n=0,1,\ldots, . 
\label{fourheun}
\end{align}
The exponents at $z=0,1$ and the position of the extra singularity are specialized by the parameter $h$, and the exponents at $z=\infty $ (i.e. $-n-3$ and $n+g+h-3$) depend on the discrete parameter $n$.
The accessory parameter ($q$ in equation (\ref{Heun})) is specialized so that the extra singularity is apparent and the equation has a polynomial solution. As shown in \ref{sec:Ia}, the equation has an extra solution for $n=-2$, \eqref{Iaextra}.

It is interesting to note that the {\em accessory parameter} $q$
in \eqref{fourheun} 
has a very special dependence on $\alpha$, $\beta$ (the exponents at $\infty$) and the location $t$ of the extra singularity:
\begin{align*}
&\frac{3(4 h-9 )n^2+3 (4 h^2-6 h-9 )n+2 h (10 h-27)}{4(h-3) (4 h-9)}\eqdef q=\alpha \beta t +q_r,\\
\alpha&\eqdef -(n+3),\quad \beta\eqdef n+g+h-3,\quad
t\eqdef -\frac{3}{4(h-3)},\quad q_r\eqdef -\frac{2h-9}{2(h-3)}.
\end{align*}
Here, the reduced accessory parameter $q_r$ does not depend on $n$.
This simplifies \eqref{fourheun} to 
\begin{align*}
 & \frac{d^2f_n}{dz^2}+\left(\frac{g+1/2}{z}+\frac{h-3/2}{z-1 }-\frac{4}{z-t} \right) \frac{df_n}{dz}  
+\left(\frac{\alpha\beta}{z(z-1)}-\frac{q_r}{z(z-1)(z-t)}\right)f_n=0.
\end{align*}
For all the examples in section \ref{sec:four}, we have 
\begin{equation}
q=\alpha\beta t+q_r,\qquad q_r:\text{does not depend on } n,
\label{qalbet}
\end{equation}
and the following Proposition.
\begin{prop} \label{heun8}
The polynomial solutions of the eight examples in section {\rm{\ref{sec:four}}} can be written in the following special form of Heun's equation with $n$-independent $q_r$:
\begin{align}
 & \frac{d^2f_n}{dz^2}+\left(\frac{\gamma}{z}+\frac{\delta}{z-1 }-\frac{4}{z-t} \right) \frac{df_n}{dz}  
+\left(\frac{\alpha\beta}{z(z-1)}-\frac{q_r}{z(z-1)(z-t)}\right)f_n=0,
\label{fourheun3}\\[2pt]
& q_r=\left\{
\begin{array}{ll}
(2+\bar{g}+\bar{h}-4\bar{g}\bar{h})/\{2(\bar{h}-1)\}& \text{Groups I
\& II},\\
(-6+\bar{g}(3-4\bar{h})+7\bar{h})/2\bar{h} & \text{Groups III
\& IV},
\end{array}\right.
\quad \gamma+\delta=\alpha+\beta+5.
\end{align}
\end{prop}
This can be verified in each case. The above relation \eqref{qalbet}
can be rephrased that the pole residues of the extra singularities are $n$-independent. This is also the case for the $X_\ell$ Jacobi polynomials which
have $\rho=-1$. See (3.11)-(3.12) in \cite{hos}.
The same statement applies to the generic multi-indexed Jacobi polynomials, (4.6)-(4.7) in \cite{os25}.

Another example is the case 
$[\tilde{\phi }^{\rm I}_{\rm 1} , \tilde{\phi } ^{\rm II}_{\rm 1} ]$
with $g=1-h$. The denominator polynomial is simply
\begin{equation*}
\mathcal{D}(\eta;1-h,h)=\eta^3\times\text{constant}.
\end{equation*}
The degree of the two-indexed Jacobi polynomial 
$\mathcal{P}_n (\eta;1-h,h )  $ is $n+2$, and the function $ (1-\eta )^{(1-h)/2} (1+\eta )^{h/2} \mathcal{P}_n (\eta;1-h,h )/\eta ^2 $ is a solution of the deformed Schr\"odinger equation.
In this case, however, the extra singularity at $\eta =0$ is located in the physical domain, and it destroys the self-adjointness of the Hamiltonian as well as the quantum mechanical interpretation.
From a viewpoint of mathematics, this case is out of the Sturm-Liouville theory, but it satisfies a second-order differential equation.
Namely, the two-indexed Jacobi polynomial $f_n=\mathcal{P}_n (\eta;1-h,h )$ satisfies
\begin{equation*}
 \frac{d^2f_n}{d\eta ^2}+\left( \frac{-h+3/2}{\eta-1}+\frac{h-3/2}{\eta+1 }-\frac{4}{\eta} \right) \frac{df_n}{d\eta } + \frac{(-n-2) (n-1) \eta + 4h-2 }{ \eta (\eta-1) (\eta+1) } f_n=0 ,
\end{equation*}
and it also satisfies Heun's differential equation by setting $\eta=1-2z$.

Similar phenomena happen in the cases $\{ [ \tilde{\phi } ^{\rm I}_{\rm 1} , \tilde{\phi } ^{\rm I}_{\rm 3} ], g=1-h \}$, $\{ [ \tilde{\phi } ^{\rm II}_{\rm 1} , \tilde{\phi } ^{\rm II}_{\rm 3} ], g=1-h \}$, $\{ [ \tilde{\phi } ^{\rm III}_{\rm 1} , \tilde{\phi } ^{\rm III}_{\rm 3}], g=h \}$.

Note that the Darboux-Crum transformation was already applied to Heun's 
equation in \cite{Tak5}. However the motivation and the results in the 
present paper are different from those in \cite{Tak5}.

\subsection{Limits to Laguerre systems}
\label{laglim}

We now discuss the relations of these Jacobi systems with the Laguerre systems presented in \cite{st}.
It is well-known that by a simple limiting procedure, one can obtain the (deformed)-radial oscillator potential 
from the (deformed)  trigonometric PT potential \cite{os19}.
For this one makes use of the limit formulas of the base polynomials
\begin{equation}
  \lim_{\beta\to\infty}P_n^{(\alpha,\,\pm\beta)}
  \bigl(1-2x\beta^{-1}\bigr)
  =L_n^{(\alpha)}(\pm x).
  \label{JtoL}
\end{equation}
By taking  the limit of infinite coupling
$h\to\infty$ with the rescaling of the coordinate:
\begin{equation}
  x=\frac{{x}^\text{L}}{\sqrt{h}},\quad
  0<x<\frac{\pi}{2}\Leftrightarrow 0<{x}^\text{L}<\frac{\pi}{2}\sqrt{h}\, ,
  \label{rescale}
\end{equation}
we have, with $\eta(x)=\cos(2x)$ and ${\eta}^\text{L}({x}^\text{L})={{x}^\text{L}}^2$, the relation
\begin{equation}
\eta (x)=1-\frac{2{\eta}^\text{L}}{h}+O(h^{-2}).
\end{equation}

Taking the limit $h\to \pm\infty$ (or equivalently, ${\bar h} \to \pm\infty$), where the upper (lower) sign applies to  Groups I and II (III and IV), we have
\begin{equation}
\lim_{h\to\pm\infty}h^{-1}\mathcal{E}_n^{\text{J}}(g,h)
  =\mathcal{E}_n^{\text{L}}(g)=4n,
\end{equation}
and
\begin{align}
\lim_{\bar{h}\to\pm\infty}\left[\frac{U^{(\mathcal{N})}(x)}{\bar{h}}\right]={\eta}^\text{L}-\frac{3}{16{\eta}^\text{L}}
+\frac{48}{4{\eta}^\text{L}\pm3}\mp\frac{288}{(4{\eta}^\text{L}\pm3)^2}-2\left(\bar{g}+\gamma+\frac12 \right),
\end{align}
with $\gamma=2$ for Groups I to III, and $\gamma=4$ for Group IV.
These are the potentials of the Laguerre systems presented in \cite{st}, where $\bar{g}=\bar{g}(\bar{h})$ takes on the respective
value for each group at the limits $\bar{h}\to \pm\infty$, namely
\begin{equation}
\bar{g}\to \frac14\; \;  ({\rm II,\  III}),\quad \frac34 \; \;  ({\rm I,\ IV}).
\label{g-limit}
\end{equation}
The Group I corresponds to the case (A), II to (B) and III to (D) in \cite{st}.

The equations (\ref{EoY12}) and (\ref{EoY34}) for $y(\eta)$ reduce to 
\begin{align}
\eta y''+(\bar{g}+\frac12-\eta)y'+\left(n + \frac{\gamma}{2} - \frac{12(4\eta\mp 3)}{(3\pm 4\eta)^2} \right)\, y=0,\quad &n=0,1, \ldots,  
\end{align}
for Groups I and II (upper sign and $\gamma=2$)  and groups III and IV (lower sign and $\gamma=2$ and $4$, respectively).
This corresponds to the equation (3.40), (3.49) and (3.55) in \cite{st}.

In the same limit, the polynomial $\mathcal{P}_n$ reduce to the corresponding  polynomial $\mathcal{L}_n$, which satisfies, from the limits of 
 (\ref{EoP12}) and (\ref{EoP34}) for $y(\eta)$:
\begin{align}
 4\eta(4\eta+3)\mathcal{L}_n^{\prime\prime}(\eta)&-2\left(-3 - 6 \bar{g} + 34 \eta - 8 \bar{g} \eta+ 8 \eta^2\right)\mathcal{L}_n^{\prime}(\eta)\n\\
 &+4\left(11 - 8 \bar{g} + 3 n + 4 (3 + n) \eta\right) \mathcal{L}_n(\eta)=0, 
\end{align}
 for Groups I and II, and
\begin{align}
  4\eta(4\eta-3)\mathcal{L}_n^{\prime\prime}(\eta)&-2\left( 3 - 2 \bar{g} (4 \eta -3) + 22 \eta + 8 \eta^2\right) \mathcal{L}_n^{\prime}(\eta)\n\\
  &+ \left( 32 - 32 \bar{g} - 6 \gamma - 12 n + 8 (4 + \gamma + 2 n) \eta \right) \mathcal{L}_n(\eta)=0,
\end{align}
for Groups III ($\gamma=2$)  and IV ($\gamma=4$).
These correspond to the equation (3.41), (3.50) and (3.56) in \cite{st}.

The above results indicate  that
the systems described in Group I, II and III  ( the potential, eigenfunctions, eigenvalues, etc.) reduce to the
 Cases (A), (B) and (D) Laguerre systems given in \cite{st}, respectively
 (for details of the Laguerre systems, see \cite{st}).  Group IV corresponds to a new Laguerre system. Since this is not 
 discussed in \cite{st},  we shall give some details
  of this case in the Appendix. 

Specifically, the forms of $U^{(\mathcal{N})}(x)$ in \eqref{potI} and \eqref{potIII} imply that the eigenfunctions have the form:
\begin{align}
\psi^{\text{J}}_n(x;g,h)=\frac{(1-\eta)^{\bar{g}/2}(1+\eta)^{\bar{h}/2}}{w_i(\eta; \bar{h})}\,\mathcal{P}_n(\eta;g,h),~~i=1,2.
\end{align}
 Under the above limiting process,
 the prefactor of $\mathcal{P}_n(\eta;g,h)$ 
(times the factor $h^{-\frac{\bar{g}}{2}} $) tends to
\begin{equation}
\frac{{\eta}^{\text{L}\frac{\bar{g}}{2}}e^{{-\frac{{\eta}^\text{L}}{2}}}}{(4{\eta}^\text{L}\pm3)^2},
\end{equation}
with $\bar{g}=3/4$ for Groups I and IV, and $\bar{g}=1/4$ for Groups II and III.
The polynomial part
 $h^{-m} \mathcal{P}_n(\eta;g,h)$ ($m=2$ for Groups I, II, and III, and $m=3$ for Group IV) approaches
 to
\begin{equation}
\mathcal{L}_n^{\text{A,\ B,\ D, \text{New}}}({\eta}^\text{L}; g),
\end{equation}
with $g=3/4, 1/4, 9/4$ and $ 13/4$, respectively.

%
%

\section{Other examples of higher degrees of apparent singularities}
\label{others}
\setcounter{equation}{0}

\subsection{The case $[ \tilde{\phi } ^{\rm II}_{\rm 2} , \tilde{\phi } ^{\rm III}_{\rm 2} ]$ and related cases}
\label{II2III2}

In the case $[ \tilde{\phi } ^{\rm II}_{\rm 2} , \tilde{\phi } ^{\rm III}_{\rm 2} ]$, the Wronskian is written as $ {\rm W} [\tilde{\phi } ^{\rm II}_{2 }, \tilde{\phi } ^{\rm III}_{2 }  ] (x)  =(1-\eta )^{3/2-g}  \mathcal{D}(\eta) $, where $\mathcal{D}(\eta )$ is a polynomial of degree $4$.
The discriminant of the polynomial $\mathcal{D}(\eta )$ is written as
\begin{align*}
-256 (2 h+1) (2 h-3) (2 h+3)^2 (2 h-5)^2 (2 g -3) (2 g -11)^3 (2 g-5)^4 \\
 \times (g+h -4) (g-h-3) (2 g h^2-2 g h-3 g-5 h^2+5 h+12)^2. 
\end{align*}
We set $ g= (5 h^2-5 h -12)/(2 h^2 -2 h-3) $, which is a root of $2 g h^2-2 g h-3 g-5 h^2+5 h+12 =0$.
Then the polynomial $\mathcal{D}(\eta )$ is written as
\begin{equation}
\{ (h+2) (h-3) \eta -( h^2 -h  +3 ) \} \{ h(h-1) \eta -(h^2 - h -3 )\} ^3 ,
\end{equation}
and the potential is 
\begin{align}
\frac{2(g-2)(g-3)}{1-\eta} + \frac{2h(h-1)}{\eta+1} -(g +h)^2 
+ \frac{24h(h-1)\{(h^2-h-3)\eta- h(h-1)\}}{ (h(h-1)\eta -(h^2-h -3 ) )^2} \n\\
+\frac{8 (h+2) (h-3) \{ (h^2-h+3) \eta - (h+2) (h-3)}{((h+2) (h-3) \eta -( h^2 -h  +3 ) )^2 }.  
\end{align}
Note that $\eta = 1 -3/\{h(h-1) \} $ is an apparent singularity of exponent $-2$ and $\eta =1 +9 /\{(h+2) (h-3) \}$ is that of exponent $-1$.

The four cases  $[ \tilde{\phi } ^{\rm II}_{\rm 2} , \tilde{\phi } ^{\rm III}_{\rm 2} ]$,  $[ \tilde{\phi } ^{\rm I}_{\rm 1} ,  \tilde{\phi } ^{\rm I}_{\rm 2} , \tilde{\phi } ^{\rm III}_{\rm 2} ]$,  $[ \tilde{\phi } ^{\rm I}_{\rm 2} ,  \tilde{\phi } ^{\rm III}_{\rm 1} , \tilde{\phi } ^{\rm III}_{\rm 2} ]$,  $[ \tilde{\phi } ^{\rm II}_{\rm 1} ,  \tilde{\phi } ^{\rm II}_{\rm 2} ,  \tilde{\phi } ^{\rm III}_{\rm 1}, \tilde{\phi } ^{\rm III}_{\rm 2} ]$
can be treated similarly to Groups I and II in section \ref{sec:fourgroups} and
the four cases  $[ \tilde{\phi } ^{\rm I}_{\rm 2} , \tilde{\phi } ^{\rm III}_{\rm 2} ]$,  $[ \tilde{\phi } ^{\rm II}_{\rm 1} ,  \tilde{\phi } ^{\rm II}_{\rm 2} , \tilde{\phi } ^{\rm III}_{\rm 2} ]$,  $[ \tilde{\phi } ^{\rm II}_{\rm 2} ,  \tilde{\phi } ^{\rm III}_{\rm 1} , \tilde{\phi } ^{\rm III}_{\rm 2} ]$,  $[ \tilde{\phi } ^{\rm I}_{\rm 1} ,  \tilde{\phi } ^{\rm I}_{\rm 2} ,  \tilde{\phi } ^{\rm III}_{\rm 1}, \tilde{\phi } ^{\rm III}_{\rm 2} ]$
corresponds to Groups III and IV.

\subsection{The case $[ \tilde{\phi }^{\rm I}_{\rm 1} , \tilde{\phi } ^{\rm II}_{\rm 2} ]$ and related cases}

In the case $[ \tilde{\phi } ^{\rm I}_{\rm 1} , \tilde{\phi } ^{\rm II}_{\rm 2} ]$, the Wronskian $  {\rm W} [\tilde{\phi } ^{\rm I}_{1 }, \tilde{\phi } ^{\rm II}_{2 }  ] (x)  = \mathcal{D}(\eta) $ is a polynomial of degree $4$.
The discriminant of the polynomial $\mathcal{D}(\eta )$ is written as
\begin{align*}
-1024 (2 h+1) (2 h-3) (2 h+3)^2 (2 g+1) (2 g-3) (2 g-5 )^2 (g-h-3)  (g-h -1)^3 \\
 \times(4 g^2 h -4 g h^2 -30 g^2-56 g h-26 h^2 +69 g +43 h -12)^2. 
\end{align*}
The equation $ 4 g^2 h -4 g h^2  -30 g^2-56 g h-26 h^2 +69 g +43 h -12 =0$ is quadratic in both $g$ and $h$.
We find some lower lying rational solutions, 
$(g, h) = (131/2,25/2)$, $(115/2, 16)$, $(113/2,41/2)$.
For these points in the parameter plane,  the polynomial $\mathcal{D}(\eta )$ is proportional to
\begin{equation}
(13+7\eta)(5\eta+7)^3,\quad (45\eta+107)(3 \eta +5)^3 ,
\quad (\eta+3)(\eta+2)^3,
\end{equation}
respectively.
The potential for the $(g,h)=(113/2,41/2)$ case is 
\begin{equation}
V(\eta)\eqdef\frac{113 \cdot 111 }{ 2(1-\eta) } +\frac{41 \cdot 39}{2(1+\eta)}-\left( \frac{113}{2} +\frac{41}2 \right) ^2 + \frac{24(2\eta +1)}{(\eta +2)^2} + \frac{8(3\eta +1)}{(\eta +3)^2}.
\end{equation}
Note that $\eta =-2 $ is an apparent singularity of exponent $-2$, $\eta =-3$ is that of exponent $-1$, and they are outside of $[-1,1]$.
The corresponding (fine-tuned) multi-indexed  polynomials 
$\{\mathcal{P}_n(\eta)\}$, 
satisfy
\begin{align}
&-(1-\eta^2)(\eta+2)(\eta+3)\frac{d^2\mathcal{P}_n(\eta)}{d\eta^2}+2 (116 + 327 \eta + 205 \eta^2 + 36 \eta^3) \frac{d\mathcal{P}_n(\eta)}{d\eta}\n\\
&\qquad +\left(-2 (146 + 355 \eta + 111 \eta^2)-4n(n+77)(\eta+2)(\eta+3)\right)\mathcal{P}_n(\eta)=0,
\end{align} 
and the solutions are
\begin{align}
\mathcal{P}_n(\eta)&=-82\left\{\left[9 (1517 + 1952 \eta + 825 \eta^2
 + 114 \eta^3)+n (n + 77) (2 + \eta)^2 (3 + \eta)
 \right]P_n^{(56,20)}(\eta)\right.\n\\
 &\qquad \qquad \left. +(1 - \eta^2) (2 + \eta) (11 + 4 \eta)\partial_\eta P_n^{(56,20)}(\eta)\right\}.
\end{align}
They satisfy the orthonormality relations
\begin{align}
&\int_{-1}^1\frac{(1-\eta)^{56}(1+\eta)^{20}}{(\eta+2)^4(\eta+3)^2}
\mathcal{P}_n(\eta)\mathcal{P}_m(\eta)d\eta=H_n\delta_{n\,m},\\
&\qquad\qquad \qquad H_n\eqdef 2^4\frac{(n(n+77) + 1102) (n (n + 77) + 1242)(n+1)_{20}(n+57)_{20}}{(2n+77)}.
\end{align}

The other parameter choices $(g,h)=(131/2,25/2),(115/2,16)$ show equally interesting characteristics. 
The cases $[ \tilde{\phi } ^{\rm I}_{\rm 3} , \tilde{\phi } ^{\rm III}_{\rm 1} ]$, $[ \tilde{\phi } ^{\rm I}_{\rm 1} , \tilde{\phi } ^{\rm III}_{\rm 3} ]$, $[ \tilde{\phi } ^{\rm II}_{\rm 3} , \tilde{\phi } ^{\rm III}_{\rm 1} ]$, $[ \tilde{\phi } ^{\rm II}_{\rm 1} , \tilde{\phi } ^{\rm III}_{\rm 3} ]$ have similar features.
%
\section{Summary and Comments}
\setcounter{equation}{0}
\label{summary}

Several families of infinitely many rational solutions are constructed of certain second order Fuchsian differential equations having four or more regular singularities.
All the extra singularities are apparent and their characteristic exponents are $-2$ or $-1$.
They are obtained by fine-tuning the parameters of multi-indexed 
Jacobi polynomials \cite{os25}.
In other words, these Fuchsian equations are derived from the Schr\"odinger equation with the P\"oschl-Teller potential through  Darboux-Crum transformations with the seed solutions having a cubic zero.
The fine-tuning of the parameters $(g,h)$, is realised as a 
curve $g=g(h)$ (or $h=h(g)$) obtained from a particular root of the discriminant of the seed solutions corresponding to a cubic zero.
These families of infinite solutions constitute orthogonal polynomials over $(-1,1)$ with the weight function of the form
$(1-x)^\alpha(1+x)^\beta/\{(ax+b)^4q(x)^2\}$, in which $q(x)$ is a polynomial in $x$.
Being the special cases of the multi-indexed Jacobi polynomials, these polynomials have ``gaps" in the degrees and they do not satisfy the three term recurrence relations.
The eight examples presented in section four are divided into four groups, each consisting of two members.
Reflecting the structure of the theory, the orthogonal polynomials of one member of the group, can be obtained from those of the other's by a simple shift of the parameters. 
Examples IIb and IIIa were discussed in \cite{gomez4} from a
different point of view.
Properties of Heun's equations with an apparent singularity are
discussed in some detail for two explicit examples.
By the standard limiting procedures leading from the Jacobi to the Laguerre
polynomials, these special cases of multi-indexed Jacobi
polynomials are shown to  reduce to the special cases of multi-indexed Laguerre polynomials. Three of them had been  explicitly constructed in \cite{st}.

It is a good challenge to construct  solutions with higher degrees ($\rho=-3,-4$) of apparent singularities within the multi-indexed Jacobi polynomials. In the preparation of 
this article, we tried in vain some combinations of lower degree
virtual state solutions.
With increasing complexity in dealing with higher degrees of
virtual state solutions, some new ideas might be needed.

\bigskip

Here is a summary table of the contents of \S4 and Appendix. The first column denotes the seed functions for the four groups introduced in \S4.
For example $ \I_2\,\III_1$ means $\text{W}[\tilde{\phi}^\I_2,\tilde{\phi}^\III_1](x)$.
The second column shows the fine tuning curve $g=g(h)$, which are linear fractional functions. The third and fourth columns are the shifted parameters $\bar{g}$ and $\bar{h}$. The fifth column gives the  equation numbers of the common forms of the denominator polynomial $\mathcal{D}(\eta)$ and the deformed potential $U^{(\mathcal{N})}(x)$.
The sixth column gives the equation numbers of the common forms of the differential equations for  the rational eigenfunction $y_n=\mathcal{P}_n/w$ \eqref{yndef} and the multi-indexed polynomial 
$\mathcal{P}_n$. Here $w$ is the reduced denominator polynomial \eqref{wdef}.
The seventh to tenth columns give the  equation numbers of the explicit forms of the 
multi-indexed polynomial 
$\mathcal{P}_n$,  the additional modes, the differential equations of the 
multi-indexed polynomial 
$\mathcal{P}_n$ and the orthogonality relations,  respectively, for the upper member of each group.

\begin{table}
\begin{center}
\medskip
Summary Table

\begin{tabular}{|c|c|c|c|c|c|c|c|c|c|}
\hline
seeds & $g(h)$ &  $\bar{g}$ & $\bar{h}$ & \ns$\mathcal{D}$\&pot.\ns & eqs. & $\mathcal{P}_n$ & add.  & 
eq.$\mathcal{P}_n$ & ort.  \\ \hline\hline
\ns $
\begin{array}{c}
 \I_2\,\III_1    \\[2pt]
 \II_1\,\II_2\,\III_1
\end{array}
$
\ns
&
\ns$
\begin{array}{c}
\frac{3(h-3)}{4h-9}     \\[4pt]
\frac{3(5h+3)}{4h+3}
\end{array}
$\ns
&
\ns$
\begin{array}{c}
 g     \\
 g-3 
\end{array}
$\ns
&
\ns$
\begin{array}{c}
 h-2    \\
 h+1  
\end{array}
$\ns
&
$
\begin{array}{l}
 (4.1) \&     \\
  (4.2)
\end{array}
$
&
\ns$
\begin{array}{c}
 y_n: (4.9)   \\
\mathcal{P}_n: (4.11) 
\end{array}
$\ns
&
\ns$
\begin{array}{c}
(\text{A}.1)      \\
\phantom{(A.4)}    
\end{array}
$\ns
&
\ns$
\begin{array}{l}
(\text{A}.2)\&    \\
(\text{A}.3) 
\end{array}
$\ns
&
\ns$
\begin{array}{c}
(\text{A}.4)      \\
\phantom{(A.4)}  
\end{array}
$\ns
&
\ns$
\begin{array}{c}
(\text{A}.5)    \\ 
\phantom{(A.4)}  
\end{array}
$\ns
\\ \hline
\ns$
\begin{array}{c}
 \I_1\,\III_2    \\[2pt]
 \II_1\,\III_2
\end{array}
$
\ns
&
\ns$
\begin{array}{c}
\frac{h}{4h-9}     \\[4pt]
\frac{9h}{4h-1}
\end{array}
$\ns
&
\ns$
\begin{array}{c}
 g     \\
 g-2  
\end{array}
$\ns
&
\ns$
\begin{array}{c}
 h-2     \\
 h 
\end{array}
$\ns
&
$
\begin{array}{l}
 (4.1) \&     \\
  (4.2)
\end{array}
$
&
\ns$
\begin{array}{c}
 y_n: (4.9)   \\
\mathcal{P}_n: (4.11) 
\end{array}
$\ns
&
\ns$
\begin{array}{c}
(\text{A}.6)      \\
\phantom{(A.4)}    
\end{array}
$\ns
&
\ns$
\begin{array}{c}
(\text{A}.7)   \\
\phantom{(A.4)}  
\end{array}
$\ns
&
\ns$
\begin{array}{c}
(\text{A}.8)      \\
\phantom{(A.4)}    
\end{array}
$\ns
&
\ns$
\begin{array}{c}
(\text{A}.9)    \\ 
\phantom{(A.4)}  
\end{array}
$\ns
\\ \hline
\ns$
\begin{array}{c}
 \II_2\,\III_1    \\[2pt]
 \I_1\,\I_2\,\III_1
\end{array}
$
\ns
&
\ns$
\begin{array}{c}
\frac{9(h-1)}{4h-3}     \\[4pt]
-\frac{3(h-3)}{4h-15}
\end{array}
$\ns
&
\ns$
\begin{array}{c}
 g-2     \\
 g+1  
\end{array}
$\ns
&
\ns$
\begin{array}{c}
 h     \\
 h-3  
\end{array}
$\ns
&
$
\begin{array}{l}
 (4.3) \&     \\
  (4.4)
\end{array}
$
&
\ns$
\begin{array}{c}
 y_n: (4.10)   \\
\mathcal{P}_n: (4.12) 
\end{array}
$\ns
&
\ns$
\begin{array}{c}
(\text{A}.11)      \\
\phantom{(A.4)}    
\end{array}
$\ns
&
\ns$
\begin{array}{l}
(\text{A}.12) \&  \\
(\text{A}.13)
\end{array}
$\ns
&
\ns$
\begin{array}{c}
(\text{A}.14)      \\
\phantom{(A.4)}   
\end{array}
$\ns
&
\ns$
\begin{array}{c}
(\text{A}.15)    \\ 
\phantom{(A.4)}  
\end{array}
$\ns
\\ \hline
\ns$
\begin{array}{c}
 \II_1\,\III_1\,\III_2   \\[2pt]
 \I_1\,\III_1\,\III_2
\end{array}
$
\ns
&
\ns$
\begin{array}{c}
\frac{3(5h-8)}{4h-7}     \\[4pt]
\frac{7h-24}{4h-15}
\end{array}
$\ns
&
\ns$
\begin{array}{c}
 g-3     \\
 g-1 
\end{array}
$\ns
&
\ns$
\begin{array}{c}
 h-1    \\
 h-3  
\end{array}
$\ns
&
$
\begin{array}{l}
 (4.3) \&     \\
  (4.4)
\end{array}
$
&
\ns$
\begin{array}{c}
 y_n: (4.10)   \\
\mathcal{P}_n: (4.12) 
\end{array}
$\ns
&
\ns$
\begin{array}{c}
(\text{A}.16)      \\
 \phantom{(A.4)}   
\end{array}
$\ns
&
\ns$
\begin{array}{l}
(\text{A}.17) -  \\
(\text{A}.19)
\end{array}
$\ns
&
\ns$
\begin{array}{c}
(\text{A}.20)      \\
\phantom{(A.4)}    
\end{array}
$\ns
&
\ns$
\begin{array}{c}
(\text{A}.21)    \\ 
\phantom{(A.4)}  
\end{array}
$\ns
\\ \hline
\end{tabular}

\end{center}
\end{table}

\section*{Acknowledgements}

C-L.\,H.  is supported in part by the
National Science Council (NSC) of the Republic of China under
Grants NSC-99-2112-M-032-002-MY3 and NSC-101-2918-I-032-001, and in part by the National Center for Theoretical Sciences (North Branch) of R.O.C.
R.\,S. is supported in part by Grant-in-Aid for Scientific Research
from the Ministry of Education, Culture, Sports, Science and Technology
(MEXT), No.23540303 and No.22540186.
K.\,T. is supported in part by the Grant-in-Aid for Young Scientists from the Japan Society for the Promotion of Science (JSPS), No.22740107.
 C-L.\,H.  would like to
thank the staff and members of YITP for the hospitality extended
to him during his visit. 

\bigskip


\appendix{\Large{\textbf{Appendix}}}

%
%
%
\section{Basic data for four cases in section 
\ref{sec:fourgroups}}
\label{datafour}


\subsection{Case Ia}
\label{sec:Ia}

\noindent $\bullet$ $g=\frac{3(h-3)}{4h-9}$,\ $\bar{g}=g$,\ $\bar{h}=h-2$.

\noindent $\bullet$ Explicit form of wave functions:
\begin{equation}
\psi^{\text{Ia}}_n (\eta;h) = \frac{ (1 - \eta)^{\frac{g}{2}}
 (1 + \eta)^{\frac{h}{2}-1}}{ (2h \eta - 6 \eta-2h +3)^2}\, \mathcal{P}_n^{\text{Ia}}(\eta; h),
   \qquad n=0,1,2,\ldots\n
\end{equation}
\begin{align}
\mathcal{P}_n^{\text{Ia}}(\eta; h)=& 
  \left(\phantom{\text{\LARGE I}}%
\!\!
1053 - 612 h + 96 h^2 - 248 h^3 + 176 h^4 - 32 h^5 - 162 n  + 
    108 h n + 144 h^2 n\right.\n\\
& - 144 h^3 n + 32 h^4 n - 162 n^2
+ 288 h n^2 - 168 h^2 n^2  + 
    32 h^3 n^2\n\\
& + (1944 - 1026 h - 1740 h^2 + 1976 h^3 - 736 h^4 +  96 h^5 + 486n - 108 h n\n\\
& \quad \ - 432 h^2 n + 240 h^3 n - 32 h^4 n +  486 n^2 - 648 h n^2 + 264 h^2 n^2 - 32 h^3 n^2) \eta \n\\
&+ (1458  - 
       4914 h + 6132 h^2 - 3520 h^3 + 944 h^4 - 96 h^5 - 216 h n \n\\
& \quad \ -  72 h^2 n + 144 h^3 n - 32 h^4 n - 216 h n^2 + 168 h^2 n^2- 32 h^3 n^2) \eta^2 \n\\
&  + (-1620 + 4248 h - 4020 h^2 + 1792 h^3-  384 h^4 + 32 h^5 - 648 n + 504 h^2 n \n\\
& \left. \quad - 240 h^3 n   + 32 h^4 n - 
       648 n^2 + 720 h n^2 - 264 h^2 n^2 + 32 h^3 n^2) \eta^3\!\!\!
       \phantom{\text{\LARGE I}}\right) P^{(g - 1/2, h - 1/2)}_n(\eta)\n\\
& 
       + (-9 + 4 h) (-1 + \eta^2) (3 - 2 h - 6 \eta + 2 h \eta) \n\\
&\quad\ \times (27 - 
    4 h^2 + 54 \eta - 30 h \eta + 4 h^2 \eta) \partial_\eta P^{(g - 1/2, h - 1/2)}_n(\eta).
\end{align}
 
 From (\ref{h-range}), one has 
 $h>9/4$ and   $h\neq 3$. The square integrability \eqref{squarecond} requires a stronger condition $h>27/10$.
 
\noindent $\bullet$ Additional mode corresponds to $n=-2$ with
\begin{align}
P^{\text{Ia}}_{-2}(\eta; h)&= 27 - 4 h^2 + 2(h-3)(2h-9)\eta,
\label{Iaextra}\\ 
\mathcal{E}^{\text{Ia}}_{-2}&=-8 (9 - 14 h + 4 h^2)/(4 h-9)=
\tilde{\mathcal{E}}_\text{1}^\aitr(\tfrac{3(h-3)}{(4h-9)},h).
\end{align}

\noindent $\bullet$ Equation satisfied by $\mathcal{P}^{\text{Ia}}_n$:
\begin{align}
&(-9 + 4 h) (3 - 2 h + 2 (-3 + h) \eta) (1 - \eta^2)y^{\prime\prime}(\eta)\n\\
& - \left(\phantom{\text{\large I}}\!\!27 + 60 h - 58 h^2 + 12 h^3 - 2 (-81 + 102 h - 48 h^2 + 8 h^3) \eta \right. \n\\
&\left. + 
   2 (-54 + 57 h - 19 h^2 + 2 h^3) (2 \eta^2 - 1)
   \phantom{\text{\large I}}\!\!\right)y^{\prime}(\eta)\n\\
&+
\left(\phantom{\text{\large I}}\!\!\!-8 h^3 (3 + n) - 8 h^2 (-20 - 3 n + n^2) - 27 (-12 + n + n^2) \right.\n\\
&\left. +      6 h (-63 + 5 n^2) + 
     2 (-3 + h) (3 + n) (18 - 18 h + 4 h^2 - 9 n + 4 h n) \eta
     \phantom{\text{\large I}}\!\!\right) y(\eta)=0.
     \label{eopB1}
\end{align}
\noindent $\bullet$ Orthonormality relations:
\begin{align}
&\int^1_{-1}\,d\eta\; \frac{(1 - \eta)^{g - 
  1/2} (1 + \eta)^{h - 5/2}} { (2h \eta - 6 \eta-2h +3)^4}  \mathcal{P}_n^{\text{Ia}}(\eta; h)\; \mathcal{P}_m^{\text{Ia}}(\eta; h)=h_n^{\text{Ia}}\delta_{nm},
  ~~~n,m=-2,0,1,2,\ldots, \n\\[2pt]
  &~~~~~ h_n^{\text{Ia}}\eqdef 2^{g+h-1}(4h-9)^2 \left(\mathcal{E}_n-\tilde{\mathcal{E}}_\text{2}^\ai\right)
\left(\mathcal{E}_n-\tilde{\mathcal{E}}_\text{1}^\aitr\right)h_n,~~~ n\neq -2.
\end{align}


\subsection{Case \text{IIa}}

$\bullet$  $g=\frac{h}{4h-9}, \bar{g}=g, \bar{h}=h-2$.

\noindent $\bullet$ Explicit form of wave functions:
\begin{align}
\psi^{\text{IIa}}_n (\eta;h)&=\frac{ (1 - \eta)^{\frac{g}{2}}
 (1 + \eta)^{\frac{h}{2}-1}}{ (2h \eta - 6 \eta-2h +3)^2}\, \mathcal{P}_n^{\text{IIa}}(\eta;g), \quad 
     n=0,1,2,\ldots,\n 
\end{align}
\begin{align}
  \mathcal{P}_n^{\text{IIa}}(\eta;h)&= 
 \left\{1458 - 918 h - 144 h^3 + 160 h^4 - 32 h^5 - 144 h n + 264 h^2 n
 - 160 h^3 n\right. \n\\
& \quad + 32 h^4 n - 162 n^2 + 288 h n^2 - 168 h^2 n^2 + 32 h^3 n^2
 + (2187 - 2106 h \n\\
&
\quad \ - 936 h^2 + 1776 h^3 - 720 h^4 + 96 h^5 + 432 h n - 600 h^2 n + 256 h^3 n\n\\
&
\quad \ - 32 h^4 n + 486 n^2- 
     648 h n^2 + 264 h^2 n^2 - 32 h^3 n^2) \eta  \n\\
& \qquad
+ (-3888 h + 6120 h^2 - 
     3624 h^3 + 960 h^4 - 96 h^5 - 192 h^2 n + 160 h^3 n\n\\
& \quad \ 
 - 32 h^4 n
- 216 h n^2 + 168 h^2 n^2 - 32 h^3 n^2) \eta^2 + (-2916 + 6048 h  \n\\
&\quad \ - 
     4932 h^2 + 1992 h^3 - 400 h^4 + 32 h^5 - 576 h n + 672 h^2 n
     - 256 h^3 n\n\\
& 
 \left. \quad \   + 32 h^4 n - 648 n^2 + 720 h n^2 - 264 h^2 n^2     + 
     32 h^3 n^2) \eta^3\right\} P^{(g - 1/2, h - 1/2)}_n(\eta)\n\\
& \quad \ 
       + (-9 + 4 h) (-1 + \eta^2)  (3 - 2 h - 6 \eta + 2 h \eta) \n\\
& \qquad \qquad \qquad \quad \times (27 - 
    4 h^2 + 54 \eta - 30 h \eta + 4 h^2 \eta) \partial_\eta P^{(g - 1/2, h - 1/2)}_n(\eta).
\end{align}
 From (\ref{h-range}), one has $h>9/4$ and  $h\neq 3$.
 The square integrability \eqref{squarecond} is satisfied.
  
\noindent $\bullet$ An additional mode corresponds to $n=-3$
with  
\begin{align}
P^{\text{IIa}}_{-3}=1,\quad {\mathcal{E}}^{\ait a}_{-3}=-\frac{12 (27 - 20 h + 4 h^2)}{-9 + 4 h}=\tilde{\mathcal{E}}^\aitr_{2}(\tfrac{h}{(4h-9)},h).
\end{align}

\noindent $\bullet$ Equation satisfied by $\mathcal{P}^{\text{IIa}}_n(\eta;h)$:
\begin{align}
&  (-9 + 4 h) (3 - 2 h + 2 (-3 + h) \eta) (1 - \eta^2)y^{\prime\prime}(\eta)\n\\
& - ((135 - 198 h + 96 h^2 - 16 h^3) \eta \n\\
& + (-3 + h) (-9 - 20 h + 
      12 h^2 + (45 - 28 h + 4 h^2) (2 \eta^2 - 1))) y^\prime (\eta)\n\\
& + (3 + n) (27 - 20 h + 4 h^2 - 9 n + 4 h n) (3 - 2 h - 6 \eta + 2 h \eta)y(\eta)=0.
\label{eopB2}
\end{align}
\noindent $\bullet$ Orthonormality relations:
\begin{align}
&\int^1_{-1}\,d\eta\; \frac{(1 - \eta)^{g - 
  1/2} (1 + \eta)^{h - 5/2}} { (2h \eta - 6 \eta-2h +3)^4}  \mathcal{P}_n^{\text{IIa}}(\eta; h)\; \mathcal{P}_m^{\text{IIa}}(\eta; h)=h_n^{\text{IIa}}\delta_{nm},
  \quad n,m=-3,0,1,2,\ldots \n\\[2pt]
& \qquad \qquad  h_n^{\text{IIa}}\eqdef 2^{g+h-1} (4h-9)^2\left(\mathcal{E}_n-\tilde{\mathcal{E}}_\text{1}^\ai\right)
\left(\mathcal{E}_n-\tilde{\mathcal{E}}_\text{2}^\aitr\right)h_{n},  ~~~n\neq -3.
\end{align}


\subsection{Case \text{IIIa}}

 \noindent $\bullet$ $g=\frac{9(h-1)}{4h-3}, \bar{g}=g-2, \bar{h}=h$.
 
 \noindent $\bullet$ Explicit form of wave functions:

 \begin{equation}
\psi^{\text{IIIa}}_n (\eta;h)=\frac{(1 - \eta)^{\frac{g}{2}-1}
 (1 + \eta)^{\frac{h}{2}}}{ (2h \eta -2h +3)^2}\, \mathcal{P}_n^{\text{IIIa}}(\eta; h),\quad n=0,1,2,\ldots, .
\end{equation} 
  \begin{align}
  \mathcal{P}_n^{\text{IIIa}}(\eta; h)&= 
\left\{405 - 648 h + 72 h^2 + 72 h^3 + 80 h^4 + 32 h^5 - 486 n + 1620 h n- 1728 h^2 n \right. \n\\
& \quad \ 
 + 432 h^3 n + 288 h^4 n - 128 h^5 n - 162 n^2 
 + 648 h n^2 - 936 h^2 n^2 + 576 h^3 n^2 \n\\
 & \quad \ 
 - 128 h^4 n^2 + (810 - 1566 h + 612 h^2 + 408 h^3 - 192 h^4 
 - 96 h^5 + 486 n - 2268 h n\n\\
 &  \quad \ 
  + 3456 h^2 n - 1584 h^3 n
   - 480 h^4 n + 384 h^5 n + 162 n^2 - 864 h n^2 + 1656 h^2 n^2 \n\\
&  \quad \
- 1344 h^3 n^2 + 384 h^4 n^2) \eta + (270 h + 36 h^2 - 528 h^3 
+ 144 h^4 + 96 h^5 + 648 h n \n\\
& \quad \
- 1944 h^2 n + 1584 h^3 n + 
     96 h^4 n - 384 h^5 n + 216 h n^2 - 792 h^2 n^2 + 960 h^3 n^2 \n\\
& \quad \
     - 384 h^4 n^2) \eta^2 + (36 h^2 + 48 h^3 - 32 h^4 - 32 h^5
     + 216 h^2 n - 432 h^3 n + 96 h^4 n\n\\
& \quad \ \left.
 + 128 h^5 n + 72 h^2 n^2
- 192 h^3 n^2 + 128 h^4 n^2) \eta^3\right\}P^{(g - 1/2, h - 1/2)}_n(\eta)\n\\
& 
       -9 (-3 + 4 h) (-1 + \eta^2)  (3 - 2 h + 2 h \eta)
         (3 - 4 h^2 - 2 h \eta + 
   4 h^2 \eta)\partial_\eta P^{(g - 1/2, h - 1/2)}_n(\eta).
 \end{align} 
 From (\ref{h-range}), one has $h<3/4$ and $h\neq 0$.
 The square integrability \eqref{squarecond} requires $-1/2<h$.
 
\noindent $\bullet$  Additional mode corresponds to $n=-2$ with
\begin{align}
 P^{\text{IIIa}}_{-2}(\eta; h)&=3 - 4 h^2 + 2 h (-1 + 2 h)\eta,\\
 \mathcal{E}^{\text{IIIa}}_{-2}&=-8 (-3 - 2 h + 4 h^2)/(-3 + 4 h)
 =\tilde{\mathcal{E}}^\aitr_1(\tfrac{3(h-1)}{(4h-9)},h). 
\end{align}

\noindent $\bullet$ Equation satisfied by $\mathcal{P}_n^{\text{IIIa}}$:
\begin{align} 
& (-3 + 4 h) (3 - 2 h + 2 h \eta) (1 - \eta^2) y^{\prime\prime}(\eta)\n\\
&+\left\{9 + 2 h^2 - 12 h^3 + 2 (9 - 6 h - 8 h^2 + 8 h^3) \eta - 
   2 h (3 - 7 h + 2 h^2) (2 \eta^2 - 1)\right\}y^\prime (\eta)\n\\
& +\left\{-8 h^3 (3 + n) - 8 h^2 (-3 + n^2) - 9 (2 + 3 n + n^2) + 
  6 h (1 + 6 n + 3 n^2)\right.\n\\
  & \qquad \qquad \qquad \qquad \qquad \left.
   + 2 h (3 + n) (4 h^2 - 3 n + h (-6 + 4 n)) \eta\right\} y(\eta)=0.
   \label{eopB3}
\end{align} 
\noindent $\bullet$  Orthonormality relations:
\begin{align} 
&\int^1_{-1}\,d\eta\; \frac{(1 - \eta)^{g - 
  3/2} (1 + \eta)^{h - 1/2}} { (2h \eta -2h +3)^4}  \mathcal{P}_n^{\text{IIIa}}(\eta; h)\; \mathcal{P}_m^{\text{IIIa}}(\eta; h)=h_n^{\text{IIIa}}\delta_{nm},
  \quad n,m=-2,0,1,2,\ldots\n\\
& \qquad \qquad \qquad \qquad  h_n^{\text{IIIa}}\eqdef 2^{g+h-1} (4h-3)^4\left(\mathcal{E}_n-\tilde{\mathcal{E}}_\text{2}^\ait\right)
\left(\mathcal{E}_n-\tilde{\mathcal{E}}_\text{1}^\aitr\right)h_{n}, ~~n\neq -2.
\end{align}

 
 \subsection{ Case \text{IVa}}

\noindent $\bullet$  $g=\frac{3(5h-8)}{4h-7}$, $\bar{g}=g-3$, $\bar{h}=h-1$.

\noindent $\bullet$ Explicit form of wave functions: 
\begin{align*}
\psi^{\text{IVa}}_n (\eta;h)&=\frac{(1 - \eta)^{\frac{g-3}{2}}
 (1 + \eta)^{\frac{h-1}{2}}}{ (2h\eta - 2 \eta - 2 h + 5)^2}\, \mathcal{P}_n^{\text{IVa}}(\eta;g),\qquad
   n=0,1,2,\ldots.   
\end{align*}
\begin{align}
  \mathcal{P}_n^{\text{IVa}}(\eta;h)=&
 \left\{133407 - 291870 h + 189396 h^2 + 30456 h^3 - 98352 h^4 + 51552 h^5- 12096 h^6 \right.\n\\
& 
+ 185640 n - 460616 h n + 440564 h^2 n - 191368 h^3 n + 24144 h^4 n + 10336 h^5 n    \n\\
&
- 4352 h^6 n
+ 54145 n^2 - 147238 h n^2+ 165212 h^2 n^2 - 98712 h^3 n^2 + 33568 h^4 n^2 \n\\
&- 6272 h^5 n^2 + 
  512 h^6 n^2 + 128 h^7 (9 + 4 n)
   + [-194400 + 574560 h - 581760 h^2 \n\\
& 
    + 141696 h^3 + 160128 h^4 - 137088 h^5 + 41472 h^6 - 164976 n + 604208 h n \n\\
&
- 830840 h^2 n + 526768 h^3 n - 129120 h^4 n - 16384 h^5 n+ 14336 h^6 n - 48118 n^2\n\\
& 
 + 187684 h n^2 - 288488 h^2 n^2 + 227232 h^3 n^2 - 97792 h^4 n^2 + 22016 h^5 n^2\n\\
&- 2048 h^6 n^2 - 
     512 h^7 (9 + 4 n)] \eta 
       + [-97686 + 4860 h + 300672 h^2 -  256608 h^3\n\\
& - 38016 h^4 + 124416 h^5 - 51840 h^6 - 67032 n - 77448 h n + 449940 h^2 n\n\\
&  - 490632 h^3 n + 194688 h^4 - 2496 h^5 n - 16896 h^6 n - 19551 n^2 - 17934 h n^2\n\\
& 
+ 134904 h^2 n^2 - 176112 h^3 n^2 + 102336 h^4 n^2 - 28416 h^5 n^2 + 3072 h^6 n^2  \n\\
&
     +768 h^7 (9 + 4 n)] \eta^2   + [32400 - 23760 h - 102240 h^2 + 
     147744 h^3   - 33408 h^4\n\\
&    - 43776 h^5 + 27648 h^6 + 49728 n- 71872 h n - 67232 h^2 n + 179392 h^3 n \n\\
&
 - 108960 h^4 n
+12800 h^5 n + 8192 h^6 n + 14504 n^2 - 24416 h n^2- 13352 h^2 n^2 \n\\
&
 + 55008 h^3 n^2 - 45568 h^4 n^2 + 15872 h^5 n^2 - 
   2048 h^6 n^2 - 512 h^7 (9 + 4 n)] \eta^3 \n\\
&     
+ [-3240 + 4320 h + 
     9576 h^2 - 21168 h^3 + 9648 h^4 + 4896 h^5 - 5184 h^6\n\\
&
 - 3360 n
+5728 h n + 7568 h^2 n - 24160 h^3 n + 19248 h^4 n  - 4256 h^5 n\n\\
&   
- 1280 h^6 n - 980 n^2 + 1904 h n^2+ 1724 h^2 n^2 -7416 h^3 n^2 + 7456 h^4 n^2 \n\\
&  \left.   - 3200 h^5 n^2 + 512 h^6 n^2 + 
     128 h^7 (9 + 4 n)] \eta^4\right\}
      P^{(g - 1/2, h - 1/2)}_n(\eta)\n\\
&       +(-7 + 4 h) (-1 + \eta^2) \left\{336 n (5 - 2 \eta)^2 (\eta-1) +98 n^2 (5 - 2 \eta)^2 (\eta-1)   \right. \n\\
&
 + 
    16 h^4 [120 
    + n (38 - 14 \eta) + 8 n^2 (\eta-1) -  69 \eta] (\eta-1)^2 
    \n\\
&
+ 32 h^5 (9 + 4 n) (\eta-1)^3 + 3 (-2048 + 4769 \eta - 2056 \eta^2 + 388 \eta^3) -8 h^3 (\eta-1)\n\\
& \quad \times[-438 + 468 \eta - 93 \eta^2 + 
       36 n (2 - 7 \eta + 4 \eta^2)+  8 n^2 (17 - 28 \eta + 11 \eta^2)] \n\\
&+ 4 h^2 [9 (-32 + 178 \eta - 152 \eta^2 + 51 \eta^3) + 6 n^2 (-143 + 333 \eta - 249 \eta^2 + 59 \eta^3) \n\\
&+ 2 n (-953 + 2373 \eta - 1896 \eta^2 + 476 \eta^3)]
- 2 h [28 n^2 (-85 + 174 \eta - 111 \eta^2 + 22 \eta^3)\n\\
& + 8 n (-895 + 1863 \eta - 1212 \eta^2 + 244 \eta^3) \n\\
&\left. + 3 (-1447 + 3859 \eta - 2144 \eta^2 + 488 \eta^3)]\right\} \partial_\eta P^{(g - 1/2, h - 1/2)}_n(\eta).
\end{align}
 From (\ref{h-range}), one has $h<7/4$ and  $h\neq 1$.
 The square integrability \eqref{squarecond} requires a stronger condition $1/2<h<13/10$.

\noindent $\bullet$ There are two  additional modes corresponding to $n=-3$ and $n=-2$:

(i) $n=-3$:   
With  
\begin{align}
P^{\text{IVa}}_{-3}= (2 h-5)\eta -2h +2,\quad
\mathcal{E}^{\text{IVa}}_{-3}=-\frac{12 (-3 + 2 h) (1 + 2 h)}{4 h-7}
=\tilde{\mathcal{E}}^{\aitr}_{2}(\tfrac{3(5h-8)}{(4h-7)},h).
\end{align}
\quad\ (ii) $n=-2$:  
With  
\begin{align}
P^{\text{IVa}}_{-2} &=107 - 88 h + 24 h^2 - 32 h^3 + 16 h^4 \n\\
&\quad - 
8 (-25 + 35 h - 14 h^3 + 4 h^4) \eta + 
4 (-1 + h)^2 (11 - 12 h + 4 h^2) \eta^2)
,\\
\mathcal{E}^{\text{IVa}}_{-2}&=-\frac{16 (-5 + 2 h^2)}{4 h-7}=
\tilde{\mathcal{E}}^{\aitr}_{1}(\tfrac{3(5h-8)}{(4h-7)},h)
\end{align}

\noindent $\bullet$ Equation satisfied by $\mathcal{P}^{\text{IVa}}_n(\eta;h)$:
\begin{align}
&  (-7 + 4 h) (5 - 2 h + 2 (-1 + h) \eta) (1 - \eta^2)y^{\prime\prime}(\eta)\n\\
&+ \left\{(-5 + 62 h - 64 h^2 + 16 h^3) \eta - (-1 + h) (19 - 28 h + 
     12 h^2 + (5 - 2 h)^2 (2 \eta^2 - 1))\right\} y^\prime (\eta)\n\\
&
+\left\{-(-5 + 2 h) [-12 - 24 n - 7 n^2 + 4 h^2 (4 + n) + 
     4 h (-4 + 2 n + n^2)]\right.\n\\
   &\quad \ \left.   
  \qquad \qquad \qquad \qquad + 2 (-1 + h) (4 + 4 h^2 + 4 h (-2 + n) - 7 n) (4 + n) \eta\right\}y(\eta)=0.
  \label{eopB4}
\end{align}

\noindent $\bullet$   Orthonormality relations:
\begin{align}
&\int^1_{-1}\!\!d\eta\frac{(1 - \eta)^{g - 
  7/2} (1 + \eta)^{h - 3/2}} { (2h\eta - 2 \eta - 2 h + 5 )^4}  \mathcal{P}_n^{\text{\text{IVa}}}(\eta; h)\; \mathcal{P}_m^{\text{\text{IVa}}}(\eta; h)=h_n^{\text{\text{IVa}}}\delta_{nm},
  ~~n,m=-3, -2, 0,1,2,\ldots,\n\\[2pt]
& ~~~~~ h_n^{\text{\text{IVa}}}\eqdef 2^{g+h-3} (4h-7)^6\left(\mathcal{E}_n-\tilde{\mathcal{E}}_\text{1}^\ait\right)
\left(\mathcal{E}_n-\tilde{\mathcal{E}}_\text{1}^\aitr\right)
\left(\mathcal{E}_n-\tilde{\mathcal{E}}_\text{2}^\aitr\right)h_n,~~~~~ n\neq -3, -2.
\end{align}

 \subsubsection{Laguerre limit}
 
 In the $h\to -\infty$ limit as described before, the 
 polynomial $ \mathcal{P}_n^{\text{IVa}}(\eta;h)$ ($n=0,1,2,\ldots$) tends to 
\begin{align}
 \mathcal{L}_n^{\text{New}}({\eta})=&\left\{5265 + (-5616 - 2592 n) {\eta} + (11232 + 5184 n) 
 {\eta}^2 \right.\n\\
 &\left.\ + (-6912 - 
       3072 n) {\eta}^3 + (2304 + 1024 n) {\eta}^4\right\}\,L_n^{(13/4)}({\eta})\n\\
       &+ \left\{1620 {\eta} - 4 (252 + 144 n) {\eta}^2 - 4 (-816 - 384 n) {\eta}^3\right. \n\\
       &\left. \qquad \qquad \qquad \qquad - 
    4 (576 + 256 n) {\eta}^4\right\} \partial_{\eta}L_n^{(13/4)}({\eta}).
\end{align}
Here we have adopted $\eta$ instead of $\eta^\text{L}$ for simplicity of presentation.
The two additional modes approach to 
\begin{align}
P^{\text{IVa}}_{-3}&\to  -(4{\eta}+3),\\
P^{\text{IVa}}_{-2} &\to  4(16{\eta}^2 + 24 {\eta} +45).
\end{align}
They  satisfy the equation
\begin{align}
& 4 {\eta} ( 4 {\eta}-3) \mathcal{L}^{\text{New}}_n({\eta})^{\prime\prime}
 -( 4 {\eta}+3) (4 {\eta}+5)  \mathcal{L}^{\text{New}}_n({\eta}) ^\prime\n\\
&~~~~~~~~ +4 (4 n {\eta}+16{\eta}-3n-4)  \mathcal{L}^{\text{New}}_n({\eta})=0,~~ n=-3,-2,0,1,2,\ldots
\end{align}




\begin{thebibliography}{99}
%

%
%
\bibitem{cod-lev}
E.\, A.\, Coddington and N.\, Levinson, {\it Theory of ordinary differential equations}, McGraw-Hill Book Company, Inc.,
New York-Toronto-London, (1955).


\bibitem{oshima}
T. Oshima, 
``Classification of Fuchsian systems and their connection problem,"
arXiv:0811.2916 [math.CA];
%
K. Hiroe and T. Oshima, ``A classification of roots of symmetric Kac-Moody root systems and its application," http://akagi.ms.u-tokyo.ac.jp/~oshima/index.html

\bibitem{os25}
S.\,Odake and R.\,Sasaki,
``Exactly Solvable Quantum Mechanics and Infinite Families of
Multi-indexed Orthogonal Polynomials,"
Phys. Lett. {\bf B702} (2011) 164-170,
{\tt arXiv:\hspace{0pt}1105.0508[math-ph]}.
For the bibliography of various exceptional orthogonal polynomials, see the references herein.

\bibitem{st}
R. Sasaki and K. Takemura, 
``Global solutions of certain second order   Fuchsian equations with    
 a high degree of apparent singularity," 
 SIGMA {\bf 8} (2012) 085 (18pp),
 {\tt arXiv:\hspace{0pt}1207.5302 [math.CA]}.


\bibitem{gomez3}
D.\,G\'{o}mez-Ullate, N.\,Kamran and R.\,Milson,
``Two-step Darboux transformations and exceptional Laguerre polynomials,"
J. Math. Anal. Appl. 387 (2012) 410-418, 
{\tt arXiv:1103.5724\hspace{0pt}[math-ph]}.


%
%

\bibitem{infhul}
L.\,Infeld and T.\,E.\,Hull,
``The factorization method,''
Rev. Mod. Phys. {\bf 23} (1951) 21-68.


\bibitem{susyqm}
See, for example, a review:
F.\,Cooper, A.\,Khare and U.\,Sukhatme,
``Supersymmetry and quantum mechanics,''
Phys. Rep. {\bf 251} (1995) 267-385.



%
%

\bibitem{darb}
G.\,Darboux,
{\it Th\'eorie g\'en\'erale des surfaces}
vol 2  Gauthier-Villars, Paris (1888).

\bibitem{crum}
M.\,M.\,Crum,
``Associated Sturm-Liouville systems,"
Quart. J. Math. Oxford Ser. (2) {\bf 6} (1955) 121-127.
{\tt arXiv:physics/9908019}.

\bibitem{adler}
M.\,G.\,Krein,
``On continuous analogue of a formula of Christoffel from the theory
of orthogonal polynomials," (Russian)
Doklady Acad. Nauk. CCCP, {\bf 113} (1957) 970-973;
V.\,\'E.\,Adler,
``A modification of Crum's method,''
Theor. Math. Phys. {\bf 101} (1994) 1381-1386.


 %
 %
 

\bibitem{Ron}
A.\, Ronveaux (ed.), {\it Heun's differential equations\/}, Oxford University Press, Oxford (1995).

\bibitem{SL}
S.\, Slavyanov and W.\, Lay, {\it Special Functions\/}, Oxford University Press, Oxford (2000).
 
 \bibitem{takemuraJPA}
K.\, Takemura,
``Heun's equation, generalized hypergeometric function and exceptional Jacobi polynomial," 
J. Phys. {\bf A45} (2012) 085211 (14pp).
{\tt arXiv:\hspace{0pt}1106.1543[math.CA]}.

%
%

 
\bibitem{gomez4}
D.\,G\'{o}mez-Ullate, N.\,Kamran and R.\,Milson,
``A conjecture on exceptional orthogonal polynomials,"
Found. Comput. Math. In press,
{\tt arXiv:\hspace{0pt}1203.6857[math-ph]}.

\bibitem{hos}
C.-L.\,Ho, S.\,Odake and R.\,Sasaki,
``Properties of the exceptional ($X_{\ell}$) Laguerre and Jacobi
polynomials,''
SIGMA {\bf 7} (2011) 107 (24pp),
{\tt arXiv:0912.5447[math-ph]}.

\bibitem{Tak5} 
K.\, Takemura,
``The Heun equation and the Calogero-Moser-Sutherland system V: 
Generalized Darboux transformations," J. Nonlinear Math. Phys. {\bf 13} 
(2006) 584--611,  {\tt arXiv:\hspace{0pt}math/0508093}.

 \bibitem{os19}
S.\,Odake and R.\,Sasaki, ``Another set of infinitely many
exceptional ($X_{\ell}$) Laguerre polynomials,'' Phys. Lett. {\bf
B684} (2010) 173,
{\tt arXiv:0911.3442[math-ph]}.
 



\end{thebibliography}
\end{document}